\documentclass{siamart190516}

\usepackage{amssymb,amsmath}
\usepackage[english]{babel}

\usepackage{mathtools}
\DeclarePairedDelimiter{\norm}{\lVert}{\rVert}

\newcommand{\R}{\mathbb R}

\newsiamremark{remark}{Remark}

\title{Finite Element Approximation of Invariant Manifolds \\
	by the Parameterization Method
	\thanks{Submitted to the editors DATE.
		\funding{J.G. and J.D.M.J. were partially supported by 
			the Sloan Foundation Grant FIDDS-17.  J.D.M.J. was partially 
			supported by the National Science Foundation grant  DMS - 1813501.}}}
\author{Jorge Gonzalez \thanks{School of Mathematics, 
		Georgia Institute of Technology, Atlanta, GA
		(\email{jgonzalez35@gatech.edu}, \url{https://people.math.gatech.edu/~jgonzalez35/index.html})}
	\and J.D. Mireles James \thanks{Department of Mathematical Sciences, 
		Florida Atlantic University, Boca Raton, FL  
		(\email{jmirelesjames@fau.edu}, \url{http://www.math.fau.edu/people/faculty/mirelesjames.php})}
	\and Necibe Tuncer  \thanks{Department of Mathematical Sciences, 
		Florida Atlantic University, Boca Raton, FL 
		(\email{ntuncer@fau.edu}, \url{http://www.math.fau.edu/people/faculty/tuncer.php})}
}

\begin{document}
\maketitle
\begin{abstract}
We combine the parameterization method for invariant manifolds with 
the finite element method for elliptic PDEs, 
to obtain a new computational framework for high order approximation of
invariant manifolds attached to unstable equilibrium solutions of nonlinear
parabolic PDEs.
The parameterization method provides an infinitesimal invariance equation 
for the invariant manifold, which we solve via a power series ansatz. 
A power matching argument leads to a 
recursive systems of linear elliptic PDEs -- the so called
homological equations -- whose solutions are the 
power series coefficients of the parameterization. 
The homological equations are solved recursively to any desired order
using finite element approximation.  The end result is a polynomial expansion
for a chart map of the manifold, with coefficients in an appropriate finite element space.  
We implement the method for a variety of example problems having both 
polynomial and non-polynomial nonlinearities, on non-convex two 
dimensional polygonal domains (not necessary simply connected), 
for equilibrium solutions with Morse indices one and two. 
We implement a-posteriori error indicators
which provide numerical evidence in support of the claim that the manifolds
are computed accurately.  
\end{abstract}



\section{Introduction} \label{sec:intro}

The present work concerns nonlinear stability 
analysis for parabolic partial differential equations (PDEs).
In particular, we
develop high order numerical methods for approximating 
local unstable manifolds attached to equilibrium
solutions of finite Morse index (finite number of unstable 
eigenvalues counted with multiplicity)
for parabolic PDEs
formulated on spatial domains with non-trivial geometry.
We show that the Taylor coefficients of an appropriate
parameterization of the local unstable manifold solve a \textit{homological 
equation} which is strongly related to the eigenvalue problem/resolvent of 
the linearization at equilibrium.    Our main goal is to 
leverage this result in the development of efficient numerical algorithms.
We stress that, since we compute the Taylor coefficients order by order  
by directly solving the homological equations, 
our method does not require numerical integration of the parabolic PDE.

Recall that the equilibrium solutions of a parabolic PDE
are found by solving the steady state equation,  and that this equation 
usually reduces to an elliptic BVP. 
Likewise, the eigenvalue problems which determine the linear stability of an 
equilibrium solution are linear elliptic BVPs of the same kind. 
Because of this, there are dramatic differences between 
parabolic problems in the case of one spatial variable and in the case of two or more.
For problems with one spatial variable, equilibrium and eigenvalue problems 
lead to two point
BVPs for ordinary differential equations (ODEs).  Such problems
are generally amenable to spectral methods (Fourier series)
which diagonalize both differential operators and multiplication 
(in Fourier and function space respectively) and which typically have
excellent convergence properties.      
Parabolic PDEs in two or more spatial variables posed on 
domains with non-trivial geometry require fundamentally different
theoretical and numerical tools. 
Finite element analysis is invaluable in this context, and   
-- since finite element methods typically employ lower regularity 
approximation schemes -- it is often necessary to study 
a weak formulation of the BVP.

Our approach is rooted in the tradition of the qualitative 
theory of dynamical systems, and exploits the parameterization method 
of Cabr\'{e}, Fontich, and de la Llave
 \cite{cabre2003parameterization1,cabre2003parameterization2,cabre2005parameterization3}.
The idea of the parameterization method is to study an auxiliary functional equation, whose solutions
correspond to chart maps of the invariant object.
The method is used widely in the field of computational 
dynamics.  The basic mathematical setup and some additional references are
discussed in Section \ref{sec:parmMethod}.
We extend the parameterization method to parabolic PDEs on non-trivial domains, and 
illustrate it's utility by implementing numerical computations
for a number of example systems.

\begin{itemize}
\item \textbf{The Fisher Equation:} scalar reaction/diffusion equation with logistic nonlinearity.
This pedagogical example illustrates the main steps of our procedure in the easiest possible setting. 
\item \textbf{The Ricker Equation:} a modification of the Fisher equation with a more realistic exponential 
nonlinearity.  We show how non-polynomial problems are treated using ideas from 
automatic differentiation for formal power series. 
\item \textbf{A modified Kuramoto-Shivisinsky Equation:} a scalar parabolic PDEs with the 
bi-harmonic Laplacian as the leading term and lower order derivatives in the nonlinearities.  
The system is a toy model of fluid dynamics. 
\end{itemize}

For each example we derive the homological equations, and implement
numerical procedures for solving them.  In the case of a non-polynomial 
nonlinearity, the necessary formal series manipulations are simplified by 
coupling the given PDEs to auxiliary equations describing the transcendental
nonlinear terms. We provide examples of this procedure, and   
develop power series expansions for 
unstable manifolds attached to equilibria with Morse indices 1 and 2.
This provides examples of computations for one and two dimensional 
unstable manifolds.    
The Fisher and Ricker Equations are nonlinear heat equations, and 
we use piecewise linear finite elements to approximate the 
coefficients of the parameterization.  
Kurramoto-Shivisinsky is a bi-harmonic Laplacian equation, 
so that higher order elements are appropriate.  
Here we utilize the Argyris element.  
We implement a-posteriori error indicators for each of the examples, 
giving evidence that the manifolds have been computed correctly.

\begin{remark}[Invariant manifolds for 1D domains]\label{rem:1d}
We remark that  Fourier-Taylor methods for computing invariant manifolds for parabolic 
problems in one spatial dimension are treated in a number of places, for example in
\cite{parmPDEone,MR4127962,MR3541499}, and
higher dimensional problems with periodic boundary conditions (including 
Dirchlet/Neumann boundary conditions on rectangles/boxes)
can also be studied using multivariate Fourier series.
We refer to the works of 
\cite{MR1430739,MR2481548,MR3720772,MR2496834,MR3640698}
for more discussion of invariant manifolds in this context.
\end{remark}
 
The remainder of the paper is organized as follows.
In Section \ref{sec:background} we review the 
finite element method for elliptic PDEs, and the parameterization 
method for invariant manifolds on Hilbert spaces.
We also provide an elementary example of the formal series
analysis for the unstable manifold in a simple finite dimensional 
example.  In Section \ref{sec:homPDE}
we extend the parameterization method to a class of parabolic
problems.  Section \ref{sec:applications}
contains the main  calculations of the paper, as we derive the homological 
equations for the main examples.  We also 
implement the recursive solution of the homological 
equations for the main examples and report on some numerical 
results.  Some conclusions and reflections are found in 
Section \ref{sec:conclusions}.


\section{Background} \label{sec:background}
While the material in this section is standard in some circles, 
the methods of the present work combine tools from 
different fields and it is worth reviewing some basic ideas.  
Our hope is that 
some brief review will help to make the paper more self contained.
The reader familiar with these ideas may want to skip ahead to 
Section \ref{sec:homPDE}, and refer back to these sections only as needed.


\subsection{The parameterization method} \label{sec:parmMethod}
The parameterization method is a general functional analytic framework for 
studying invariant manifolds, originally developed for 
fixed points for maps on Banach spaces 
\cite{MR1976079,MR1976080,MR2177465}, and 
for whiskered tori of quasi-periodic maps
\cite{MR2289544,MR2299977,MR2240743}.
Since then it has been extended to a number of settings for both 
discrete and continuous dynamical systems, in both finite and infinite 
dimensions.  A complete overview of the literature is beyond the scope 
of the present brief introduction, and the interested reader will find a much more 
complete overview -- including a wealth of references to the literature -- 
in the recent book on the topic \cite{MR3467671}.
Several papers more closely related to the present work include 
works of \cite{MR3736145,MR3501842,MR3808251} on delay differential equations,
KAM for PDEs \cite{MR3900818}, and 
unstable manifolds for PDEs defined on compact intervals 
\cite{parmPDEone}, and on the whole line \cite{MR4127962}.
More recently the parameterization method has been used to
develop a mathematically rigorous approach to  
optimal mode selection in nonlinear 
model reduction by projecting onto 
spectral submanifolds 
\cite{MR3800257,MR3832468,MR4233869}. This research direction has been further developed and combined with large finite element systems demonstrating its potential for industrial applications \cite{vizzaccaro2021high, opreni2022fast}.

\subsubsection{Parameterization method for vector fields on Hilbert spaces} \label{sec:parmPDEs}
We give a brief review the parameterization method, in the context
of evolution problems on Hilbert spaces.  The main application we have in mind is
the dynamics of a semi-flow generated by parabolic PDE.  In particular, we discuss the invariance
equation for the local unstable manifold attached to an  
equilibrium solution.   

Let $\mathcal{H}$ be an $L^2$ Hilbert space
and 
$F \colon \mathcal{H} \to L^2$ 
be a Frechet differentiable mapping.  In fact, 
we only require that the derivative of $F$ at each point 
is densely defined.  Consider the evolution equation 
\begin{equation} \label{eq:evolution}
\frac{\partial}{\partial t} u(t) = F(u(t)),  \quad \quad \quad 
\mbox{with } u(0) \in \mathcal{H} \mbox{ given}.
\end{equation}
An orbit segment (or \textit{solution curve}) for 
Equation \eqref{eq:evolution}
 is a 
smooth curve $\gamma \colon (a, b) \to \mathcal{H}$
having
\[
\frac{d}{dt} \gamma(t) = F(\gamma(t)), 
\]
for each $t \in (a, b)$.  If $b = \infty$ then $\gamma$ is a 
said to be a full forward orbit.  Since $F$ dose not depend on time, we can 
always choose $a = 0$.

The simplest type of orbits are equilibria, that is, solutions which do 
not change in time.  
For $u_0 \in \mathcal{H}$, the curve $\gamma(t) = u_0$ 
is a constant solution of Equation \eqref{eq:evolution}
if and only if
\[
F(u_0)  = 0.
\]
For a given equilibrium solution $u_0$,  we would like to understand 
first it's linear stability, and then it's nonlinear stability.
That is, we would like to understand how orbits in a 
neighborhood of $u_0$ escape from that neighborhood.
  
Let $A = DF(u_0)$, and define the Morse index of $u_0$ to be the 
number of unstable eigenvalues of $A$, counted with multiplicity.  
We assume that Equation \eqref{eq:evolution} is parabolic, so that   
$A$ generates a compact semi-group
$e^{At}$.  This insures that the Morse index of $A$ is finite.
Let $\lambda_1, \ldots, \lambda_M$ denote the unstable eigenvalues 
ordered so that 
\[
0 < \mbox{real}\left(\lambda_1\right) \leq
\ldots \leq \mbox{real}\left(\lambda_M\right).
\]
Suppose for the sake of simplicity that each unstable eigenvalue has multiplicity 
one, and that they are all real (though both assumptions can be removed --
see \cite{MR1976079,MR3518609}), 
and let $\xi_1, \ldots, \xi_M \in \mathcal{H}$ denote associated 
eigenfunctions, so that 
\[
A \xi_j = \lambda_j \xi_j, \quad \quad \quad 1 \leq j \leq M.
\]

Suppose that $\gamma \colon (-\infty, 0] \to \mathcal{H}$ is a solution curve
for Equation \eqref{eq:evolution} and that $u \in \mathcal{H}$.
We say that $\gamma$ is an infinite pre-history
for $u$, accumulating in backward time to the 
equilibrium $u_0$, if 
\[
\gamma(0) = u, 
\quad \quad \quad \mbox{and} \quad \quad \quad
\lim_{t \to -\infty} \gamma(t) = u_0. 
\] 
The unstable manifold attached to $u_0$, denoted 
$W^u(u_0)$, is the set of all $u \in \mathcal{H}$
which have an infinite pre-history, accumulating at $u_0$.
The intersection of $W^u(u_0)$ with a neighborhood $U$ of $u_0$ 
is called a local unstable manifold for $u_0$, and is denoted by 
\[
W^s(u_0) \cap U = W_{\mbox{\tiny loc}}^u(u_0, U).
\]
By the unstable manifold theorem, there exists a neighborhood $U$ of $u_0$ so that 
$W_{\mbox{\tiny loc}}^u(u_0, U)$ is a smooth manifold, 
diffeomorphic to an $M$-disk, and 
tangent to the unstable eigenspace of $A$ at $u_0$.
Moreover, if $A$ is hyperbolic (that is, if $A$ has no eigenvalues on the 
imaginary axis), then 
$W_{\mbox{\tiny loc}}^u(u_0, U)$ is the set of all $u \in U$
which have well-defined backwards history remaining in a neighborhood 
of $u_0$ for all time $t \leq 0$.

\begin{figure}[!t]
\centering
\includegraphics[width=3.5in]{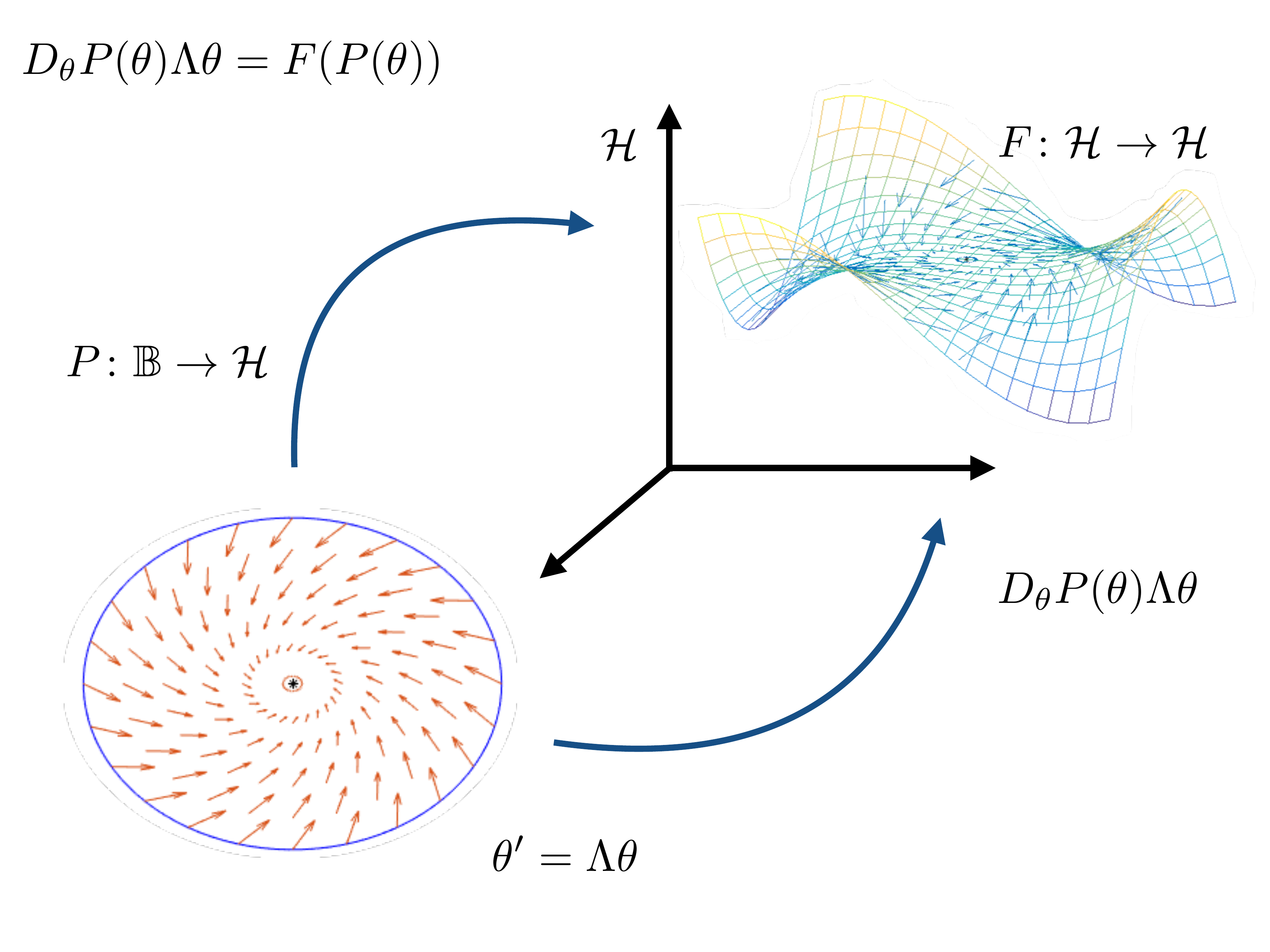}
\caption{Schematic representation of the invariance equation given in 
Equation \eqref{eq:invEq}.  The idea is the $DP$ pushes forward the 
vector field $\Lambda$ modeling the dynamics on the 
unstable manifold.  This push forward should be equal, on the 
image of $P$, to the vector field $F$ generating the 
full dynamics.}
\label{fig:parmMethodSchematic}
\end{figure}

We are now ready to introduce the parameterization method.
Let $\mathbb{B} = [-1, 1]^M$ denote the $M$-dimensional unit hypercube.  
We seek a $P \colon \mathbb{B} \to \mathcal{H}$ having that 
\begin{equation}\label{eq:constraint1}
P(0) = u_0, 
\end{equation}
\begin{equation} \label{eq:constraint2}
\partial_j P(0) = \xi_j, \quad \quad \quad \quad 1 \leq j \leq M,
\end{equation}
and that 
\[
P\left([-1,1]^M\right) \subset W^u(u_0, U),
\]
for some open set $U$ containing $u_0$.
Any such $P$ is local unstable manifold attached to $u_0$.
Since any reparameterization of $P$ is again a parameterization 
of a local unstable manifold, the problem has infinitely many freedoms and 
we need to impose an additional (infinite dimensional) constraint
to isolate a single parameterization. 

Write  
\begin{equation} \label{eq:defLambda}		
\Lambda = 
\left(\begin{array}{ccc}
\lambda_1 & \ldots & 0 \\
\vdots & \ddots & \vdots \\
0 & \ldots & \lambda_{M}
\end{array}\right).
\end{equation}
The main idea of the parameterization method is to look for $P$ which, 
in addition to satisfying  the constraint Equations \eqref{eq:constraint1} and 
\eqref{eq:constraint2}, is a solution of the \textit{invariance equation}
\begin{equation} \label{eq:invEq}
F(P(\theta))=DP(\theta)\Lambda \theta,
\quad \quad \quad \mbox{for all } \theta \in \mathbb{B} = [-1,1]^M.
\end{equation}
We remark that the choice of ``unit'' domain is a normalization which will become
more clear as we proceed.  

\begin{figure}[!t]
\centering
\includegraphics[width=3.5in]{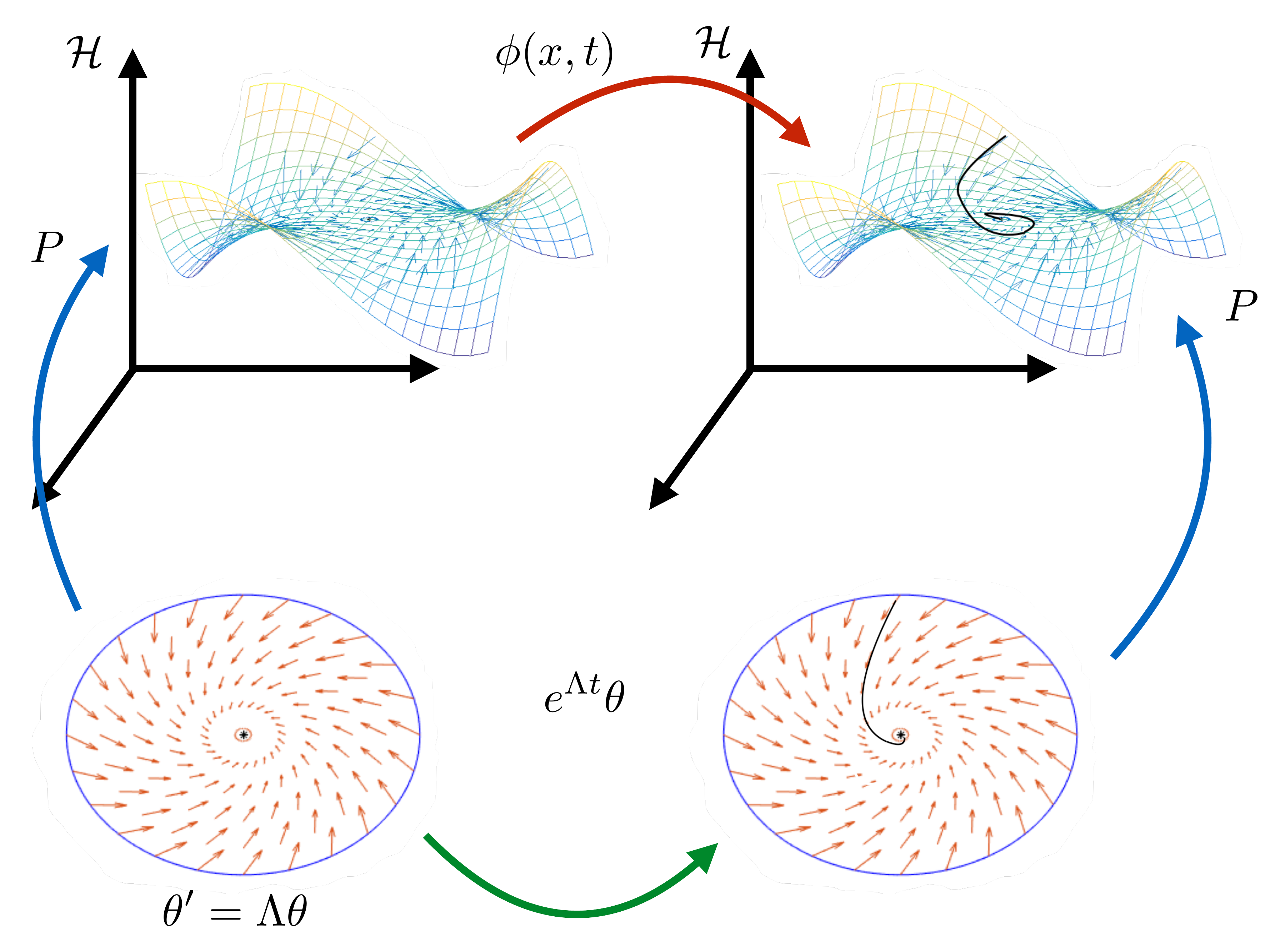}
\caption{The orbit correspondence induced by the invariance Equation.
The orbits generated by the vector field $\Lambda$ accumulate in backwards
time to the origin in $\mathbb{B}$.  Then $P$ lifts these orbits to orbits in 
$\mathcal{H}$ which accumulate at the equilibrium $u_0$.  From this it 
follows that image of $P$ is a local unstable manifold. 
\eqref{eq:invEq}}
\label{fig:infConj}
\end{figure}

Figure \ref{fig:parmMethodSchematic} illustrates the geometric meaning of
Equation \eqref{eq:invEq}.  The equation requires that the push forward of the linear
vector field $\Lambda$ by $DP$
equals the vector field $F$ restricted to the image of $P$.
Loosely speaking, since the two vector fields match on the 
image of $P$ they must generate the same dynamics --
with the dynamics generated by $\Lambda$ well understood.
We then expect that  $P$ maps orbits of $\Lambda$ in $\mathbb{B}$ 
to orbits of $F$ on the image of $P$.  
Since $P$ maps orbits to orbits, 
Equation \eqref{eq:invEq} is called an infinitesimal conjugacy 
equation.  The geometric meaning of Equation \eqref{eq:invEq} 
is illustrated in Figure 
\ref{fig:infConj}, is made precise by the following lemma.

%
%

\begin{lemma}[Orbit correspondence] \label{lem:orbitCorr}
Assume that the unstable eigenvalues $\lambda_1, \ldots, \lambda_M$
are real and distinct.
Suppose that $P \colon [-1,1]^M \to \mathcal{H}$ satisfies the first order 
constraints of Equations \eqref{eq:constraint1} and \eqref{eq:constraint2},
and that $P$ is a smooth solution of Equation \eqref{eq:invEq} on 
$\mathbb{B} = (-1,1)^M$.
Then $P$ parameterizes a local unstable manifold for $u_0$.    
\end{lemma}

\begin{proof}
First observe that since the
domain $\mathbb{B}$ is a topological disk, and $P$ is a smooth
mapping, we have that the image of $P$ is a smooth $M$-dimensional 
manifold.   Also observe that the constraint given in Equation  \eqref{eq:constraint2}
implies that $P$ is tangent to 
the unstable eigenspace of $DF(u_0)$ at $u_0$.

Now fix $\theta \in (-1, 1)^M$, and define the curve 
$\gamma_{\theta} \colon (-\infty, 0] \to \mathcal{H}$ by 
\[
\gamma_{\theta}(t) = P\left(e^{\Lambda t} \theta \right).
\]
We observe that $\gamma_{\theta}$ 
is a solution curve for $F$.  
To see this, we first note that $\gamma_\theta$ is well 
defined for all backward time, as 
for all $t \in (-\infty, 0]$ we have that 
\[
\hat \theta := e^{\Lambda t} \theta \in \mathbb{B}.
\]
This is because the entries of $\Lambda$ are unstable, real, and distinct.  
To see that $\gamma_\theta(t)$ solves the differential equation, note that 
\begin{align*}
\frac{d}{dt} \gamma_{\theta}(t) & = \frac{d}{dt} P\left(e^{\Lambda t}  \theta\right) \\
& = D P\left( e^{\Lambda t}  \theta \right) \Lambda e^{\Lambda t}\theta \\
& = D  P\left( e^{\Lambda t}  \theta \right) \Lambda e^{\Lambda t}\theta \\
& = D P (\hat \theta) \Lambda \hat \theta \\
& = F(P(\hat \theta)) \\
& = F(P (e^{\Lambda t} \theta)) \\
& = F(\gamma_{\theta}(t)), 
\end{align*}
as desired.

In addition to being a solution curve, we have that 
$\gamma_{\theta}$ accumulates at $u_0$ 
in backward time.  To see this, we simply compute the limit
\begin{align*}
\lim_{t \to - \infty} \gamma_{\theta}(t) &= \lim_{t \to - \infty} P\left(e^{\Lambda t} \theta\right) \\
&= P\left(  \lim_{t \to - \infty} e^{\Lambda t} \theta\right) \\
&=P\left(  0 \right) \\
& = u_0,
\end{align*}
where we have used the assumption that $P$ is smooth, and hence
continuous on $[-1,1]^M$.  
Since $\theta$ was arbitrary, we see that every point $P(\theta)$ on 
the image of $P$ has a backward orbit which accumulates 
at $u_0$.  That is 
\[
\mbox{image}(P) \subset W^u(u_0).
\] 
Since $\mbox{image}(P)$ is an $M$-dimensional disk containing $u_0$ and contained in 
the unstable manifold, we have that $\mbox{image}(P)$ is a local unstable manifold as desired.  
\end{proof}

We remark that if $F$ generates a semi-flow $\Phi$ near $u_0$, then Lemma \ref{lem:orbitCorr}
says that $P$ satisfies the flow conjugacy 
\begin{equation} \label{eq:conj}
P(e^{\Lambda t} \theta) = \Phi(P(\theta), t), 
\end{equation}
for all $t$ such that $e^{\Lambda t} \theta \in (-1,1)^M$.
That is, $P$ conjugates the flow generated by $\Lambda$ to 
the flow generated by $F$.

\begin{remark}[Complex conjugate unstable eigenvalues] \label{rem:complexConjEigs}
Complex conjugate eigenvalues are easily incorporated into this set-up 
by choosing associated complex conjugate eigenfunctions and proceeding as above.
This results in complex conjugate coefficients for the parameterization $P$.  The use of 
complex conjugate variables (in the appropriate components of $\theta$) results in 
$P$ having real image, i.e. recovers the parameterization of the real manifold.  The 
only difference is that one has to adjust the domain of the parameterization in 
the variables corresponding to the complex conjugate eigenvalues, choosing 
unit disks instead of unit intervals.  
In this sense the PDE case is no different from the ODE case described in detail in 
\cite{MR3207723}, where the interested reader can find more a complete discussion.
\end{remark}

\subsubsection{Formal solution of Equation \eqref{eq:invEq}: an ODE example} \label{sec:tutorial}

In this section we we illustrate the use of the
parameterization method as a computational tool
for a simple example.
The idea is to develop a formal series solution of Equation \eqref{eq:invEq}.
Such formal calculations 
play a critical role in the remainder of the present work, and
are much more involved for PDEs than for ODEs.
To separate those complications which are inherent to the 
method from those which are due to 
PDEs,  we explain the procedure for the planar 
vector field $F \colon \mathbb{R}^2 \to \mathbb{R}^2$ 
(Hilbert space is the plane) given by 
\begin{equation} \label{eq:exODE} 
F(x,y) = \left(
\begin{array}{c}
x + y \\
1 - x^2
\end{array}
\right).
\end{equation}
We are interested in the orbit structure of $\mathbb{R}^2$ generated by 
the ODE
\[
\frac{d \gamma}{dt} = F(\gamma), 
\]
where
\[
\gamma(t) = \left(
\begin{array}{c}
x(t) \\
y(t)
\end{array}
\right).
\]
Note for future use that  
\begin{equation} \label{eq:exDir}
DF(x,y) = \left(
\begin{array}{cc}
1 & 1 \\
-2x & 0
\end{array}
\right).
\end{equation}
Suppose that $p_0 \in \mathbb{R}^2$ has $F(p_0) = 0$, 
so that $p_0$ is an equilibrium solution of the ODE.  Suppose further that 
$DF(p_0)$ has one unstable eigenvalue $\lambda > 0$ and that 
the remaining eigenvalue is stable.  Let $\xi \in \mathbb{R}^2$
denote an eigenvector associated with $\lambda$.

We look for a function 
$P \colon [-1, 1] \to \mathbb{R}^2$ with 
\[
P(0) = p_0
\quad \quad \mbox{and} \quad \quad 
 P'(0) = \xi,
\]
parameterizing the one dimensional unstable manifold attached to $p_0$.  
In the one dimensional case 
the invariance equation of Equation \eqref{eq:invEq}
reduces to 
\begin{equation} \label{eq:1DinvEq}
\lambda \theta \frac{d}{d \theta} P(\theta) = F(P(\theta)),
\end{equation}
for $\theta \in (-1, 1)$.
We look for a power series solution 
of Equation \eqref{eq:1DinvEq} 
of the form 
\[
P(\theta) = 
\left(
\begin{array}{c}
P_1(\theta) \\
P_2(\theta)
\end{array}\right) = 
\sum_{n = 0}^\infty \left(
\begin{array}{c}
a_n \\
b_n 
\end{array}
\right) \theta^n,
\]
and impose first order constraints
\[
\left(
\begin{array}{c}
a_0 \\
b_0 
\end{array}
\right) = p_0, \quad \quad \mbox{and} \quad \quad 
\left(
\begin{array}{c}
a_1 \\
b_1 
\end{array}
\right) = \xi.
\]

To work out the higher order coefficients we note that, 
on the level of formal power series, the left hand side of Equation 
\eqref{eq:1DinvEq} is 
\begin{equation} \label{eq:DP_lokta}
\lambda \theta \frac{d}{d\theta} P(\theta) = 
\sum_{n = 0}^\infty \lambda n  \left(
\begin{array}{c}
a_n \\
b_n 
\end{array}
\right) \theta^n,
\end{equation}
and that the right hand side of Equation \eqref{eq:1DinvEq} is 
\begin{eqnarray}
F(P(\theta)) &=  \left(
\begin{array}{c}
P_1(\theta) + P_2(\theta) \\
1 - P_1(\theta)^2
\end{array}
\right) \nonumber \\
& = \sum_{n = 0}^\infty 
\left(
\begin{array}{c}
a_n +  b_n \\
\delta_n - \sum_{k=0}^n a_{n-k} a_k
\end{array}
\right) \theta^n. \label{eq:FP_lokta}
\end{eqnarray}
Here we have used the Cauchy product formula
for the coefficients of $P_1(\theta)^2$, and defined 
\[
\delta_n = \begin{cases}
1 & n = 0 \\
0 & n \geq 1
\end{cases}.
\]

Returning to the invariance Equation \eqref{eq:1DinvEq}, we  
set the right hand side of Equation \eqref{eq:DP_lokta} equal 
to Equation \eqref{eq:FP_lokta}, match like powers of 
$\theta$, and recall the definition of $\delta_n$  to obtain  
\begin{equation} \label{eq:hom1Ex1}
\lambda n  \left(
\begin{array}{c}
a_n \\
b_n 
\end{array}
\right) = \left(
\begin{array}{c}
a_n +  b_n \\
- \sum_{k=0}^n a_{n-k} a_k
\end{array}
\right),
\end{equation}
for $n \geq 1$.
We seek to isolate terms of order $n$, and derive a equation for $p_n$
in terms of lower order coefficients.
Since there are still some terms order $n$ locked in the sum, we note that 
for $n \geq 2$
\[
\sum_{k = 0}^n a_{n-k} a_k = 
2 a_0 a_n + \sum_{k = 1}^{n-1} a_{n-k} a_k,
\]
where the new sum on the right contains no terms of order $n$.
Exploiting this identity, Equation \eqref{eq:hom1Ex1} becomes 
\begin{equation*} 
n \lambda   \left(
\begin{array}{c}
a_n \\
b_n 
\end{array}
\right) = \left(
\begin{array}{c}
a_n + b_n \\
-2 a_0 a_n - \sum_{k = 1}^{n-1} a_{n-k} a_k
\end{array}
\right),
\end{equation*}
or 
\begin{equation*} 
\left(
\begin{array}{c}
a_n + b_n - n\lambda a_n \\
-2 a_0 a_n  - n \lambda b_n
\end{array}
\right)
= \left(
\begin{array}{c}
 0 \\
 \sum_{k = 1}^{n-1} a_{n-k} a_k
\end{array}
\right).
\end{equation*}
This is 
\[
\left[
\begin{array}{cc}
1  - n\lambda & 1 \\
-2 a_0                      & - n \lambda
\end{array}
\right] \left(
\begin{array}{c}
a_n \\
b_n
\end{array}
\right) = 
\left(
\begin{array}{c}
  0\\
 \sum_{k = 1}^{n-1} a_{n-k} a_k
\end{array}
\right),
\]
which, after referring back to Equation \eqref{eq:exDir}, we rewrite as 
\begin{equation} \label{eq:homEqEx}
\left(DF(a_0, b_0) - n\lambda \mbox{Id} \right) p_n = s_n,  \quad \quad \quad 
n \geq 2,
\end{equation}
where 
\[
p_n = \left(
\begin{array}{c}
 a_n \\
 b_n
\end{array}
\right), \quad \quad \mbox{and} \quad \quad 
s_n = \left(
\begin{array}{c}
  0 \\
 \sum_{k = 1}^{n-1} a_{n-k} a_k
\end{array}
\right).
\]
Again, note that $s_n$ depends only on terms of order less than $n$.

We refer to Equation \eqref{eq:homEqEx} as the \textit{homoloical equations} for $P$,
and note that they are linear algebraic equations for the power series coefficients of
the parameterization.  
We now ask, \textit{are the homological equations solvable?}
To answer this we note that since $P(0) = p_0 = (a_0, b_0)$ is an equilibrium solution, 
the left hand side of Equation \eqref{eq:homEqEx}
is the characteristic matrix for the derivative $DF(p_0)$.  
The characteristic matrix is invertible if and only if 
$n \lambda$ is not an eigenvalue of $DF(p_0)$.  
Since $\lambda > 0$, and since the remaining eigenvalue of $DF(p_0)$
is negative, we see that for $n \geq 2$,  $n \lambda$ is never an eigenvalue.
Then the homological equations are uniquely solvable to all orders, and 
the power series solution of Equation \eqref{eq:1DinvEq}, when $F$ is given by 
 Equation \eqref{eq:exODE}, is formally well defined.

This implies that the coefficients of $P$ are uniquely determined after the 
first order data (equilibrium and eigenvector) are fixed.
Then the only freedom in determining the solution is 
the choice of the scaling of the eigenvector $\xi$.  This non-uniqness
is used to control the growth rate of the coefficients of $P$, providing 
numerical stability.

\begin{remark}[Non-resonance and the parameterization method] \label{rem:nonRes}
The condition 
\begin{equation} \label{eq:nonResODE}
n \lambda \notin \mbox{spec} DF (p_0) \quad \quad \quad n \geq 2, 
\end{equation}
is  called a non-resonance condition.  In fact it is an \textit{inner non-resonance
condition} as we are computing the unstable manifold, and Equation \eqref{eq:nonResODE}
involves linear combinations of the (in this case unique) unstable eigenvalues. 
We will see in Section \ref{sec:homPDE} that the non-resonance 
conditions are similar, but somewhat more subtle for higher dimensional 
unstable manifolds.  
\end{remark}

\begin{remark}[Stable manifolds for ODEs] \label{rem:stableManifolds}
Note that replacing $\lambda$  with a stable eigenvalue
in the above discussion changes nothing.  This reflects the general fact that 
in finite dimensions, the parameterization method applies
equally well to both stable and unstable manifolds.  However, 
an equilibrium solution of a parabolic PDE typically has infinitely many 
stable eigenvalues which make it impossible to overcome the non-resonance 
conditions. This is why the present work focuses on unstable manifolds
for parabolic PDEs.  
\end{remark}

\begin{figure}[!t]
	\begin{center}
	\includegraphics[scale=0.5]{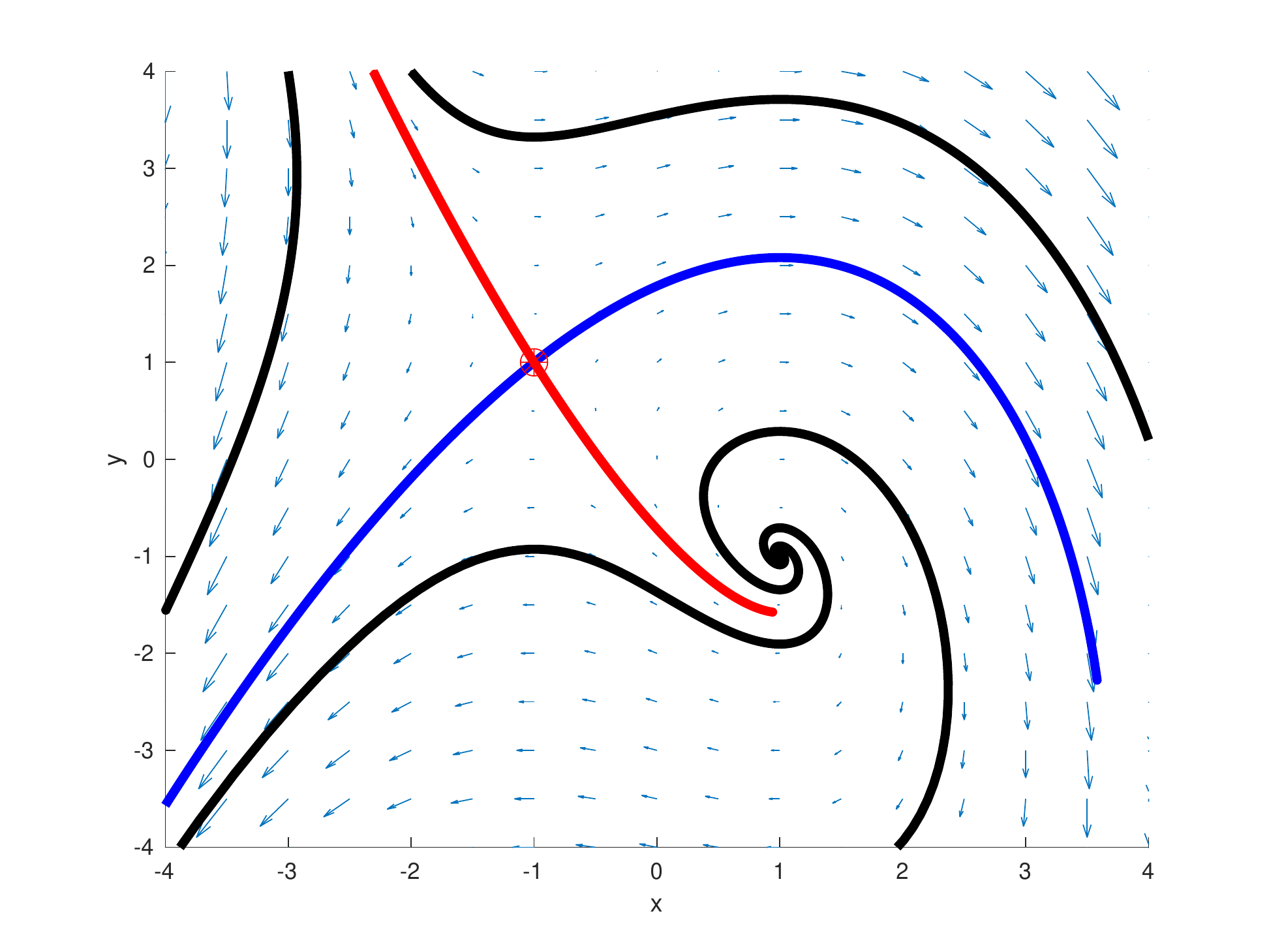}
	\end{center}
	\caption{\textbf{Stable/unstable manifold visualization:} 
	dynamics generated by the vector field given in Equation \eqref{eq:exODE}.
	Several reference orbits are illustrated by black curves.  These are obtained by
	numerical integration of several arbitrarily chosen initial conditions.  The main features
	of the phase space are the 
	saddle equilibrium at $(-1,1)$ and the repelling equilibrium at $(1,-1)$.  
	We compute the local unstable and local stable manifold parameterizations 
	$P^N$ and $Q^N$ for the saddle stable equilibrium $(-1,1)$ to order $N = 100$.
	The unstable and stable eigenvectors to lengths of $13$ and $10.5$ respectively,
	The images $P^N([-1,1])$ and $Q^N([-1,1])$ are plotted as blue (unstable) 
	and red (stable) curves.   In both cases the plots of the manifolds are generated only by 
	plotting the approximating polynomials: the manifolds are not extended using numerical 
	integration. This illustrates that it is often possible to approximate a substantial portion of 
	the unstable manifold using the parameterization method. 
	(Of course numerical integration could be used to extend the manifolds even further). 
	We observe that the unstable manifold parameterization (blue curve)
	follows a ``fold'', that is, the curve is not the graph over the unstable eigenspace
	of any function.  The stable manifold on the other hand seems have been approximated 
	up to very near it's maximal 
	radius of convergence, as computing additional terms has very little effect on the picture, 
	and we are not able to reach a fold.  
 	} \label{fig:odeExManifolds}
\end{figure}

\subsubsection{A numerical example}
The vector field of Equation \eqref{eq:exODE} has 
equilibrium solutions $f(x_{1,2}, y_{1,2}) = (0,0)$ at 
\[
\left(
\begin{array}{c}
x_1 \\
y_1
\end{array}
\right) = 
\left(
\begin{array}{c}
-1 \\
1
\end{array}
\right), \quad \quad \mbox{and} \quad \quad 
\left(
\begin{array}{c}
x_2 \\
y_2
\end{array}
\right) = 
\left(
\begin{array}{c}
1 \\
-1
\end{array}
\right),
\]
and one can check that 
\begin{equation}\label{eq:DF_ex}
DF(-1, 1) = 
\left(
\begin{array}{cc}
1 & 1 \\
2 & 0
\end{array}
\right),
\end{equation}
has eigenvalues $2, -1$.  Hence the 
equilibrium $(-1, 1)$ is a hyperbolic saddle.  
Let $\lambda = 2$ denote the 
unstable eigenvalue.  One can check that 
\[
\xi = \left(
\begin{array}{c}
1 \\
1
\end{array}
\right),
\]
is an associated unstable eigenvector.  

The zero-th and first order terms of the parameterization are
\[
\left(
\begin{array}{c}
a_0 \\
b_0
\end{array}
\right) = 
\left(
\begin{array}{c}
-1 \\
1
\end{array}
\right) \quad \quad \mbox{and} \quad \quad 
\left(
\begin{array}{c}
a_1 \\
b_1
\end{array}
\right) = 
\left(
\begin{array}{c}
1 \\
1
\end{array}
\right),
\]
and the second order term is determined by solving the homological equation 
of Equation \eqref{eq:homEqEx} with $n = 2$ as follows.  Recalling 
the definition of $s_n$, and noting that 
$a_1 = 1$, when $n = 2$ we have that 
\[
\left. \sum_{k = 1}^{n-1} a_{n-k} a_k  \right|_{n = 2} = a_1^2 = 1,
\]
and that 
\[
s_2 = 
\left(
\begin{array}{c}
0 \\
1
\end{array}
\right).
\]
Moreover, since $\lambda = 2$ and $a_0 = -1$ we 
recall Equation \eqref{eq:DF_ex}, and 
have that  
\[
DF(-1,1) - 2 \lambda \mbox{Id} = \left[
\begin{array}{cc}
1 - 2 \lambda & 1 \\
2 & - 2 \lambda
\end{array}
\right] = 
 \left[
\begin{array}{cc}
-3 & 1 \\
2 & - 4
\end{array}
\right].
\]
Solving 
\[
\left[DF(-1,1) - 2 \lambda \mbox{Id} \right] p_2 = s_2,
\]
gives
\[
p_2 = 
\left(
\begin{array}{c}
-0.1 \\
-0.3
\end{array}
\right).
\]
From this we conclude that the second order local unstable manifold approximation is  
\begin{equation} \label{eq:parmEx1}
P^2(\theta) = \left(
\begin{array}{c}
 -1 \\
 1 
\end{array}
\right) + 
\left(
\begin{array}{c}
 1 \\
 1 
\end{array}
\right) \theta  +
\left(
\begin{array}{c}
  -0.1  \\
 - 0.3 
\end{array}
\right) \theta^2.
\end{equation}
Third and higher order terms are computed recursively following 
the same recipe.  

Roughly speaking, how accurate is the approximation above?  
Since the remainder term in the approximation given by $P^2$ in Equation 
\eqref{eq:parmEx1} is cubic in $\theta$, we expect that the size of the truncation error 
has 
\[
E_2(\theta) = \|P(\theta) - P^2(\theta) \| \leq C |\theta|^3, 
\]
for some constant $C$.  Suppose that we restrict the domain of our parameterization to 
\[
\theta \in [-10^{-5}, 10^{-5}].
\]
Then $E_2$ is of order $(10^{-5})^3 = 10^{-15}$, so that the size of 
the truncation error is roughly 5 multiples machine precision.
In practice, we prefer to rescale the length of the eigenvector, and 
take the domain of $P^N$ normalized to a unit cube.  
See the following remark.

\begin{remark}[Rescaling the eigenvector to optimizing the coefficient decay]
\label{rem:rescaling}
Suppose now that we compute the coefficients of $P^N$ to order $N = 20$, 
using the same eigenvector $\xi = [1,1]$.  
Rather than listing the resulting coefficients order by order, 
we remark that the coefficients decay like 
\[
\| p_n\| \approx 65 \times 10^{-1.18 n},
\]
(found by taking an exponential best fit algorithm) 
and that
\[
\|p_{20}\| \approx  1.56 \times 10^{-22},
\]
a quantity far smaller than machine precision.  
Note that coefficients below machine precision do 
not contribute (numerically) to the approximation, and this 
is wasted effort.

To obtain a more significant result, we increase the scaling of 
the unstable eigenvector, taking $P'(0) = s \xi$
with some $s > 1$.
For example, rescaling the eigenvector by $s = 2.5$ 
and recomputing the coefficients leads to a $20$-th order 
polynomial whose final coefficient vector has magnitude $1.4 \times 10^{-14}$.  
Since the final coefficient is close to, but still above machine precision  
-- and hence numerically significant -- this choice of scaling is
nearly optimal for the order $N=20$ calculation.

Experimenting a little more in this way, we find that 
taking $s = 13$, and computing the parameterization to order $N = 100$, 
gives coefficients which decay exponentially fast and in such at way that  
the last coefficient had magnitude roughly machine epsilon.
A plot illustrating the results of the order $N = 100$ calculation 
is given in Figure \ref{fig:odeExManifolds}.
Note that the unstable manifold, which is shown as the blue curve, is
not the graph  of a function over the tangent space (span of the eigenvector).  
This illustrates the well known fact that the parameterization method
can ``follow folds'' in the manifold.  
The reader interested in more refined approaches to choosing the 
computational parameters in the parameterization method might consult 
\cite{MR3437754}, where methods for optimizing the calculations under 
certain constraints are discussed in detail.  
\end{remark}

\begin{remark}[Visualization in a Function space] \label{rem:visualization}
 The parameterization method is extremely valuable for 
 visualizing invariant manifolds when the dimension of the 
 phase space is low.  However the remainder of the paper 
 concerns infinite dimensional problems, and visualization is much more problematic.
 For the parabolic PDEs studied below, the phase space is a Sobolev space, 
 and each point on the manifold is actually a function represented as a linear 
 combination of finite elements.  In this setting 
 it is more natural to plot the points on the 
 manifolds as functions over the given domain. That is, we visualize the manifold as a curve or surface 
 of functions. Nevertheless, 
 it is valuable to keep in mind the picture in Figure \ref{fig:odeExManifolds}
when trying to interpret the results.  
\end{remark}


\subsection{Finite element methods for elliptic linear elliptic PDE} \label{sec:finiteElements}
In this section we briefly review the basics of  
finite element analysis for elliptic BVPs
needed for our numerical implementations.
Excellent reference for this now classic material include 
\cite{MR0520174,MR1930132,MR2597943}.
Let $\Omega \subset \mathbb{R}^d$ denote an open set 
and let $\mathcal{H}(\Omega)$ be an  $L^2$ Sobolev space on $\Omega$ 
(hence a Hilbert space).  Let  $\mathcal{H}^{\vee}$ denote the dual space 
consisting of all bounded linear functionals on $\mathcal{H}$.

Consider a uniformly elliptic linear PDE of the form
\[
\mathcal{L}u= f,
\]
having  boundary conditions $\{ B_i(u)|_{\partial \Omega} = g_i\}$.
We ask that $\mathcal{L}$ be a densely defined 
linear operator, $u \in \mathcal{H}(\Omega)$, and $f \in L^2(\Omega)$.
The $B_i$'s  denote \textit{boundary operators}, for example
directional derivatives, 
or more complicated constraints at the boundary, $g_i\in  L^2(\partial \Omega)$.
 
A weak formulation of the problem is obtained after multiplying the equation by a
$v\in \mathcal{H}(\Omega)$, applying 
Green's formula (integration by parts), and imposing the
boundary conditions.
This results in the variational problem 
\begin{equation} \label{eq:ellipticProblem}
\text{Find} \hspace{0.2cm} u\in \mathcal{H} \hspace{0.2cm} \text{such that} \hspace{0.2cm} \forall v \in \mathcal{H}, \hspace{0.2cm} \langle u,v \rangle_{\mathcal{L}}=\langle f, v \rangle,
\end{equation}
where 
\[
\langle f, v \rangle = \int_\Omega f v,
\]
and $\langle u,v \rangle_{\mathcal{L}}$ is a bilinear form derived 
from $\mathcal{L}$ as described above (Green's formula/boundary conditions). 
The classical \textit{Lax-Milgram lemma} insures that 
the problem has a unique solution $u$, assuming that 
$\langle \cdot,\cdot \rangle_{\mathcal{L}} \colon \mathcal{H} \times \mathcal{H} \to \mathbb{R}$
is a continuous $\mathcal{H}$-elliptic bilinear form and 
$\langle f,\cdot \rangle \colon \mathcal{H} \to \mathbb{R}$ is a bounded 
linear functional (i.e, $\langle f,\cdot \rangle \in \mathcal{H}^{\vee}$).

The finite element method (FEM) is a Galerkin projection
approach to numerically solving Equation \eqref{eq:ellipticProblem},
and consists of three main steps:
\begin{itemize}
\item[1.] Triangulate $\Omega \subset \mathbb{R}^d$:  obtain 
(often polygonal) mesh which discretizes the problem domain.  
\item[2.] Choose interpolants for $\mathcal{H}$ on the mesh:
construct a basis for the interpolant space where the basis functions
have nearly disjoint support over mesh elements.  This is the finite 
element basis and it's span is a finite element space.
\item[3.] Solve the sparse linear system obtained by projecting the 
the weak formulation of the PDE (Equation \eqref{eq:ellipticProblem})
onto the finite element basis.  This reduces the problem to numerical 
linear algebra.
\end{itemize}

\begin{figure}[t!]
\begin{center}
	\begin{minipage}[b]{0.32\textwidth}
		\includegraphics[width=\textwidth]{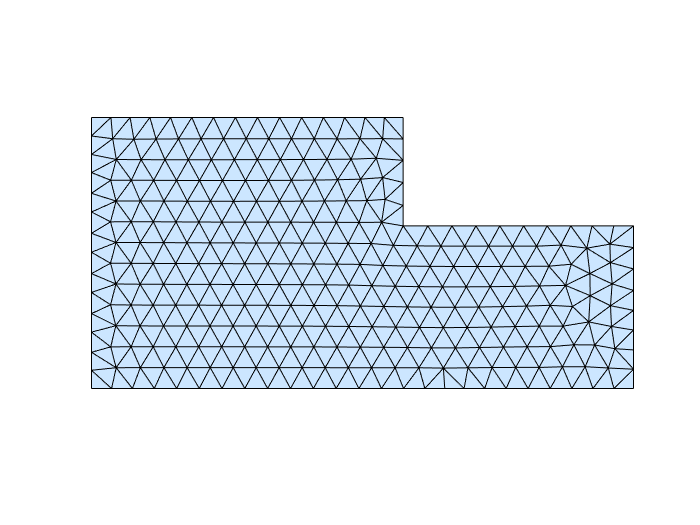}
	\end{minipage}
	\begin{minipage}[b]{0.32\textwidth}
		\includegraphics[width=\textwidth]{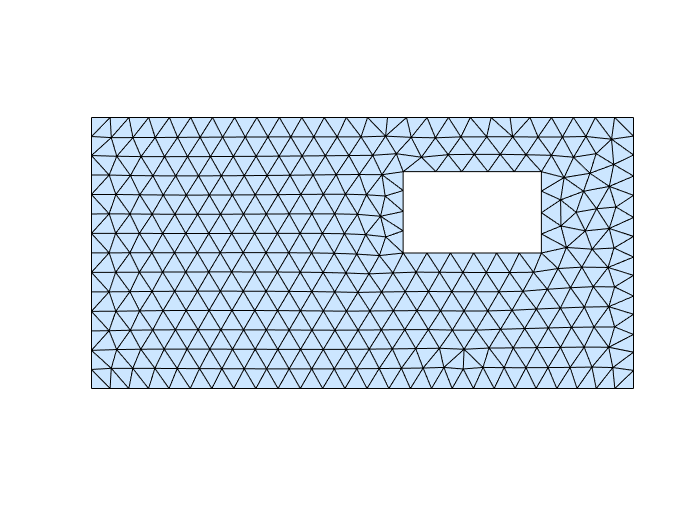}
	\end{minipage}
	\begin{minipage}[b]{0.32\textwidth}
		\includegraphics[width=\textwidth]{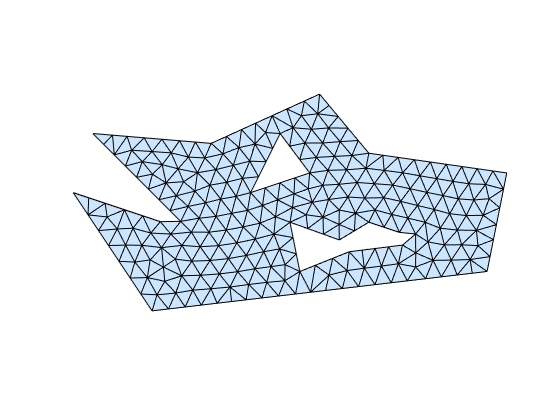}
	\end{minipage}
\end{center}
\caption{Three example domains used in this paper.  Note that they are
non-convex, and non-simply connected.
Left: the $\textit{L}$ domain: it has a reemergent corner.
Center: the \textit{Door} domain: not simply connected.  
Right: the \textit{Polygon with holes} domain: toy model of
a ``natural'' domain like a lake with islands.   
}%
\label{fig:planarDomains}%
\end{figure}

In the present work we focus on $\Omega\subset \R^2$ a polygonal domain.
However, we do not 
require $\Omega$ to be convex or even simply connected. 
More precisely, we use the  
domains illustrated in Figure \ref{fig:planarDomains}.
The next three subsections discuss the three steps above.

\subsubsection{Triangulation of $\Omega\subset \R^2$}
 
Let $\{T_i\}_{i = 1}^{{ne}}$ denote the elements of the triangulation so that  
\[
\Omega=\bigcup\limits_{i=1}^{{ne}} T_i.
\]
Here $T_i$ is the $i^{th}$ triangle, and ${ne}$ is the number of 
triangular elements.  We require that if the boundary of two triangles meet, 
then their intersection must be at a common edge. 
We remark that other discretizations can be considered, 
for example as in the Bogner-Fox-Schmit elements \cite{MR0520174} (quadrilaterals), 
or even a combination of rectangles and triangles. 
Also, the discretization does not need to be regular but can be adapted to the model 
and domain, leading to more efficient approximations.

\subsubsection{Constructing the basis elements}

The basis elements, which are required to have  
``small" compact support in $ \mathcal{H}$,  are typically 
chosen to be piecewise polynomial. In this paper we use  linear 
polynomials for second order problems (Laplacian operator),
and fifth degree polynomials for some degree 4 examples (Bi-harmonic Laplacian). 
These Argyris elements are discussed in more detail in  Section \ref{sec:KS_manifolds}. 
More general basis elements can be considered such as 
special rational functions (for example Zienkiewicz triangles \cite{MR0520174}).

A finite element is denoted by $E=[z_1, \cdots ,z_{{nn}}] \subset T$, 
where $T$ is an arbitrary triangle and the $z_i$ are \textit{control points} or \textit{nodes}.
$S_i:=\{L_{ij}: 1\leq j \leq s_i\}$ denotes a corresponding sets of \textit{control operators} 
evaluated at $z_i$ ({nn} is the number of nodes in $T$ and $s_i$ denotes the total 
number of operators assigned to the node $z_i$). Typically, the nodes consist of the 
vertices along with a few other carefully chosen points. In general, they are not required to be uniformly distributed in $T$. 

Denote by ${lnb}:=\sum\limits_{k=1}^{{nn}} s_i$ the total number of operators associated
with the element $E$. These letters appropriately stand for ``local number of basis" since the
operators are used to determine the basis elements associated with $T$. Let 
$\mathbb{B} \subset \mathcal{H}$  denotes the span of the 
basis elements and define 
\[
S=\bigcup\limits_{i=1}^{{nn}} S_i = \{L_i : 1\leq i \leq \sum\limits_{k=1}^{{nn}} s_i\}.
\]
Then for each $k$,  the system $L(\phi):= (L_1(\phi), \cdots, L_{{lnb}}(\phi))= e_k$ 
has a unique solution in $\mathbb{B}$.  Here $e_k$ is the $k^{th}$ elementary basis 
vector in $\mathbb{R}^{{lnb}}$.

Let $V_h:=\mbox{span}\{\phi_i\}_{i=1}^{{nb}}\subset  \mathcal{H}$ denote 
an interpolation space for $\mathcal{H}$, where ${nb}$ is the total number 
of basis elements.
We want that 
\[
\mathbb{B}=\mathbb{P}_k:=\{p: p  \hspace{0.2cm} \text{is a polynomial of degree at most}
 \hspace{0.2cm} k \},
 \] 
so,  must have ${lnb}=\frac{(k+1)(k+2)}{2}$. Imposing regularity conditions (for example
continunity) on the solution $u$ imposes further restrictions on the elements. 
For $\mathbb{B}=\mathbb{P}_1$ the elements are of the 
form $E=[n_1, n_2, n_3]$ where the $n_i$'s are the 
vertices of the triangles, and $S_i=\{{id}\}$ for all $i$'s, with ${id}(\phi)(n_i)=\phi(n_i)$.

\subsubsection{Computing the projection}
Let $u \in \mathcal{H}$ denote the solution of Equation \eqref{eq:ellipticProblem}.
The projection of $u$ into $V_h$ is found by 
solving a weak formulation of Equation \eqref{eq:ellipticProblem} on $V_h$. More precisely, 
write $u_h=\sum\limits_{i=1}^{{nb}} c_i \phi_i$ 
and solve the linear system
\[
\sum_{j=1}^{{nb}}c_j \langle \phi_j, \phi_i \rangle_{\mathcal{L}} = \langle f, \phi_i \rangle.
\]

It follows by an application of the Lax-Milgram lemma that the 
matrix $ \Big( \langle \phi_j, \phi_i \rangle_{\mathcal{L}} \Big)$
is invertible.

In general, a Lagrange type interpolation of a function $f$ over $T$ with control set $\{S_i\}$ 
is given by 
\[
 \Pi_T(f)=\sum\limits_{i=1}^{{lnb}}L_i (f)(z_{n(i)}) \frac{det(A_i)}{det(A)} 
\]
where $L_i\in S=\bigcup\limits_{i=1}^{{nn}} S_i$, and 
the index $n(i)=k$ for $i$ such that  $s_0+\cdots s_{k-1} + 1 \leq i \leq s_0+\cdots s_{k}$. 
Here we define 
$A_{ij}=\Big( L_i(x^m y^n)(z_{n(i)})\Big)$, 
and ${(A_k)}_{ij}= (1-\delta_{ki})A_{ij} + \delta_{ki} L_i( x^my^n)$, 
where $j=\frac{(m+n)(m+n+1)}{2} + (n+1)$. Let $S_0:=\emptyset$ for 
convenience of expressing $n(i)$.

For low order polynomial bases the integrals can be evaluate exactly. 
For higher order bases it is often more practical to use quadrature rules of 
sufficiently high degree to approximate the integrals.  Such rules have the form  
\[
\int_{\Omega} f=\sum_{i=1}^{{ne}}\int_{T_i} f \approx \sum_{i=1}^{{{ne}}} \sum_{j=i}^{{nq}} w_{j}^{T_i} f(q_{j}^{T_i}),
\]
where ${nq}$ is the degree of the quadrature rule, $q_j^{T_i}$ are the quadratures points, 
and $w_j^{T_i}$ are some appropriately chosen weights. 
Then 
\[
\langle \phi_j,\phi_i \rangle_{\mathcal{L}}^{q} c^q= \langle f, \phi_i \rangle^{q},
\]

where $\langle \cdot,\cdot \rangle_{\mathcal{L}}^{q}$ and $\langle f,\cdot \rangle^{q}$ denote
the quadrature approximation of the bilinear form and linear functional respectively. 
If $\langle \cdot,\cdot \rangle_{\mathcal{L}}$ is $ \mathcal{H}$-elliptic, 
it follows that $\langle \cdot,\cdot \rangle_{\mathcal{L}}^{q}$ is $V_h$-elliptic, which 
implies that $\Big( \langle \phi_j,\phi_i \rangle_{\mathcal{L}}^{q} \Big)$ is invertible.
In general, the $ \mathcal{H}$-elliptic property of 
$\langle \cdot,\cdot \rangle_{\mathcal{L}}$ is 
established using the Sobolev embedding theorems/Poincar\'{e} inequalities. 

For any polynomial basis there is ${nq}$ large enough so that  $\langle \phi_j,\phi_i \rangle_{\mathcal{L}}^{q}=\langle \phi_j,\phi_i \rangle_{\mathcal{L}}$,  
in which case 
\[
\norm{c^q - c} \leq \left\|
\Big (\langle f, \phi_i \rangle^{q} - \langle f, \phi_i \rangle \Big)\right\| 
\left\|\Big( \langle \phi_j,\phi_i \rangle_{\mathcal{L}}^{q} \Big)^{-1} \right\|.
\] 

Approximating $f=p+\epsilon$, with $p$ polynomial, we have 
 \[
 \left\| \Big( \langle f, \phi_i \rangle^{q}  - \langle f, \phi_i \rangle  \Big) \right\|  \leq 2 \sup(\epsilon)
\left\|\Big(\langle 1, \phi_i \rangle \Big) \right\|,
\]
for ${nq}$ large enough.

Bounding the projection error for a polynomial basis of order $k$ requires 
assumptions about the domain $\Omega$.  It follows, for example, by the 
the Bramble-Hilbert lemma that $\norm{u-u_h}_{1,\Omega} = O(h^k)$, 
where $u_h$ denotes the projection of the solution $u$ to the finite dimensional 
vector space $V_h$.  Of course, more sophisticated and practical ways of estimating 
these errors can be found in the literature.


\section{Formal power series and the homological equations for parabolic PDEs} \label{sec:homPDE}
We now turn to the main problem of this paper, which is to extend the kinds of  
calculations illustrated in Section \ref{sec:tutorial} to the ``vector fields'' 
on Sobolev spaces generated by parabolic PDEs.  To this end we 
introduce a fairly simple class of nonlinear heat equations which we find sufficient to 
highlight the main issues.  Nevertheless, the discussion in this section generalizes
to parabolic equations involving more general
elliptic operators,  to problems formulated on spatial
domains of three or more dimensions with more general boundary conditions, 
and even to systems of PDEs.  
Indeed, our goal in this section is not to describe the most general possible setting
 but rather 
to illustrate the application parameterization method, and especially the 
solution of Equation \eqref{eq:invEq}, for an interesting class of PDEs.
Some extensions are given in Section \ref{sec:applications}.

Let $\Omega \subset \mathbb{R}^2$ denote bounded, planar, polygonal domain  and 
$f \colon \mathbb{R} \times \Omega \to \mathbb{R}$
be a smooth function.  
Consider the class of scalar parabolic PDEs given by 
\begin{equation} \label{eq:reactionDiffusion}
\frac{\partial}{\partial t} u(t, x,y) = \Delta u(t, x, y) + f(u(t, x, y),  x, y), 
\end{equation}
with the Neumann boundary conditions
\[
\frac{\partial}{\partial \mathbf{n}} u(t, x, y) = 0 \quad \mbox{for } (x, y ) \in \partial \Omega.
\]
Fix $\mathcal{H} = H^1(\Omega)$. We are interested in the 
dynamics of the semi-flow generated by the vector field
 $F \colon \mathcal{H} \to L^2(\Omega)$  given 
by 
\[
F(u) = \Delta u + f(u, x,y).
\]
Note that $F$ maps a dense subset of $\mathcal{H}$ into 
$\mathcal{H}$.
 
We now consider an equilibrium solution.  That is, 
suppose that $u_0 \colon \Omega \to \mathbb{R}$ is in $\mathcal{H}$ and 
is a solution of the weak form of the elliptic nonlinear boundary 
value problem  
\[
\Delta u(x, y) + f(u(x, y),  x, y) = 0,
\]
subject to the Neumann boundary conditions. 
More precisely, this means that $u_0$ satisfies 
\[
-\int_{\Omega} \nabla u(x,y) \cdot \nabla \phi(x,y) + \int_{\Omega} f(u, x, y) \phi(x,y) = 0,
\]
for all $\phi \in \mathcal{H}$. 

Suppose also that $u_0$ has Morse index $M$.  That is, we assume that 
$\lambda_1, \ldots, \lambda_M\in (0, \infty)$ are the unstable eigenvalues,
each with multiplicity one.  Let 
$\xi_1, \ldots, \xi_M \colon \Omega \to \mathbb{R}$ denote 
associated unstable eigenfunctions, i.e. 
solutions in $\mathcal{H}$ of the weak form of the eigenvalue problem
\[
\Delta \xi(x,y) + \partial_1 f(u_0,x,y) \xi = \lambda \xi(x,y),
\] 
again subject to the boundary conditions.  

We look for $P \colon [-1, 1]^M \to \mathcal{H}$ solving 
Equation \eqref{eq:invEq}, 
with $P$ given by the formal power series 
\[
P(\theta_1, \ldots, \theta_M, x,y) = 
\sum_{n_1 = 0}^\infty \ldots \sum_{n_M = 0}^\infty p_{n_1, \ldots, n_M}(x,y) \theta_1^{n_1} \ldots \theta_M^{n_M}.
\]
Here each coefficient $p_{n_1, \ldots, n_M} \in \mathcal{H}$ 
is required to satisfy the boundary conditions.  
Moreover, imposing the constraints of Equations \eqref{eq:constraint1} and \eqref{eq:constraint2} 
gives that the first order coefficients of $P$ are
\[
p_{0, \ldots, 0}(x,y) = u_0(x,y), 
\]
and 
\[
p_{1, \ldots, 0}(x,y) = \xi_1(x,y), 
\quad \quad \quad \ldots \quad \quad \quad
p_{0, \ldots, 1}(x,y) = \xi_M(x,y).
\]

To work out the higher order coefficients we follow the blueprint of 
Section \ref{sec:tutorial}.  Begin by letting 
$\Lambda$ denote the diagonal matrix of unstable eigenvalues
as in Equation \eqref{eq:defLambda}.
Calculating the push forward of $\Lambda$ by $DP$ on the 
level of power series gives  
\begin{eqnarray*}
DP(\theta, x,y) \Lambda \theta &= \left[
\partial_1 P(\theta, x,y), \ldots, \partial_M P(\theta, x,y) \right]
\left(
\begin{array}{c}
\lambda_1 \theta_1 \\
\vdots \\
\lambda_M \theta_M 
\end{array}
\right) \\
& = 
\lambda_1 \theta_1 \frac{\partial}{\partial \theta_1} P(\theta, x,y) + 
\ldots + \lambda_M \theta_M \frac{\partial}{\partial \theta_M} P(\theta, x,y) \\
&= 
\sum_{n_1 = 0}^\infty \ldots \sum_{n_M = 0}^\infty 
(n_1 \lambda_1 + \ldots + n_M  \lambda_M)
p_{n_1, \ldots, n_M}(x,y) \theta_1^{n_1} \ldots \theta_M^{n_M}.
\end{eqnarray*}
Observe that the value of this series at $\theta = 0$ is zero.

Next  consider  
\begin{align*}
F(P(\theta, x,y)) &= \Delta P(\theta, x, y) + f(P(\theta, x,y), x, y).
\end{align*}
Formally speaking, the Laplacian commutes with the infinite sum, and we have that  
\[
\Delta P(\theta, x, y) =  \sum_{n_1 = 0}^\infty \ldots \sum_{n_M = 0}^\infty 
\Delta p_{n_1, \ldots, n_M}(x,y) \theta_1^{n_1} \ldots \theta_M^{n_M}.
\]
If $f$ is analytic then $f(P(\theta, x,y), x,y)$ 
admits a power series representation. 
(For $f$ only $C^k$ regularity the argument below is modified accordingly). 
Let us write 
\[
f(P(\theta, x,y), x,y) = \sum_{n_1 = 0}^\infty \ldots \sum_{n_M = 0}^\infty 
q_{n_1, \ldots, n_M}(x,y) \theta_1^{n_1} \ldots \theta_M^{n_M},
\]
where the $q_{n_1, \ldots, n_M}$ are the formal Taylor coefficients of the 
composition, and each depends on the coefficients of $P$.  Efficient computation 
of the $q_{n_1, \ldots, n_M}$ 
best illustrated through examples in the next section and  for the moment we
remark that, for any given multi-index $(n_1, \ldots, n_M) \in \mathbb{N}^M$, 
the dependence of $q_{n_1, \ldots, n_M}$ on $p_{n_1, \ldots, n_M}$
has 
\[
q_{n_1, \ldots, n_M} = D_1 f (u_0, x,y) p_{n_1, \ldots, n_M} + S_{n_1, \ldots, n_M},
\]
where $S_{n_1, \ldots, n_M}$ depends only on coefficients of $P$ of lower order.
This follows from the Fa\'{a} di Bruno formula.

Matching like powers in Equation \eqref{eq:invEq} leads to 
\begin{align*}
&(n_1 \lambda_1 + \ldots + n_M  \lambda_M)
p_{n_1, \ldots, n_M}\\
& =  \Delta  p_{n_1, \ldots, n_M} + q_{n_1, \ldots, n_M}\\
& = \Delta  p_{n_1, \ldots, n_M}  +
D_1 f (u_0, x,y) p_{n_1, \ldots, n_M} + S_{n_1, \ldots, n_M},
\end{align*}
so that 

\begin{gather*} 
\Delta  p_{n_1, \ldots, n_M} +
D_1 f (u_0, x,y) p_{n_1, \ldots, n_M} - (n_1 \lambda_1 + \ldots + n_M  \lambda_M)
p_{n_1, \ldots, n_M}\\
= - S_{n_1, \ldots, n_M}.
\end{gather*}
That is, $p_{n_1, \ldots, n_M}$ solves the linear equation 
\begin{equation} \label{eq:homEqHighLevel}
\left( 
DF(u_0) - (n_1 \lambda_1 + \ldots + n_M  \lambda_M) \mbox{Id}_{\mathcal{H}}
\right)
p_{n_1, \ldots, n_M} = - S_{n_1, \ldots, n_M},
\end{equation}
where the right hand side depends only on lower order terms.

Equation \eqref{eq:homEqHighLevel} is \textit{the homological equation} for 
the unstable manifold for $F$ at $u_0$.  
Observe that Equation \eqref{eq:homEqHighLevel} is a linear elliptic PDE with the same boundary 
conditions as the original reaction/diffusion equation \eqref{eq:reactionDiffusion}. Indeed, 
the linear operator on the left hand side is the resolvent of $DF(u_0)$,
evaluated at the complex numbers $n_1 \lambda_1 + \ldots + n_M  \lambda_M$.  
Then each Taylor coefficient of $P$ is the solution of a linear problem 
no more complicated than the linearized equation at $u_0$, so that these equations are themselves
amiable to finite element analysis under mild assumptions on the domain $\Omega$.  

This is a general fact which makes the parameterization method so useful.
The homological equations determining the jets of the invariant manifold parameterization
 are always linear equations in the same category as
the steady state equations for the equilibrium solution itself.  
For example when considering a finite dimensional problem in Section 
 \ref{sec:tutorial}, the steady state equations were systems of 
 $n$ nonlinear algebraic  
 equations in $n$ unknowns, and in this case the
 homological equations turned out to be systems of $n$ 
 linear equations in $n$ unknowns.
 Moreover, the homological equations involved the characteristic 
 matrix for the derivative of the vector field at the equilibrium.  
 
 In the calculations just discussed, the steady state equation is 
a nonlinear elliptic BVPs, and the homological equations turn out the 
be linear elliptic BVPs on the same domain with the same boundary conditions.
In fact the linear operator is just the resolvent of the differential, in direct 
analogy with the finite dimensional case.
In the remarks below, we expand on several similarities between the results just 
derived and the simple example calculation considered in Section \ref{sec:tutorial}.

\begin{remark}[Non-resonance conditions and existence of a formal solution]
Observe that Equation \eqref{eq:homEqHighLevel} has a unique solution 
if and only if the \textit{non-resonance condition} 
\begin{equation} \label{eq:nonRes}
n_1 \lambda_1 + \ldots + n_M  \lambda_M \notin \mbox{spec}(DF(u_0)),
\end{equation}
is satisfied whenever $n _1 + \ldots + n_M \geq 2$.
Since $\lambda_1, \ldots, \lambda_M$ are the only unstable eigenvalues of $DF(u_0)$,
and since $DF(u_0)$ generates a compact semi-group,  
we have that the countably many remaining eigenvalues are stable.  
Since the  $n_1, \ldots, n_M$ are all positive, there are only finitely many opportunities 
for $n_1 \lambda_1 + \ldots + n_M  \lambda_M$ to be an eigenvalue.  
If Equation \eqref{eq:nonRes} is satisfied for all multi-indices $(n_1, \ldots, n_M) \in \mathbb{N}^M$
with $n_1 + \ldots + n_M \geq 2$ then we say that the unstable eigenvalues are \textit{non-resonant}, 
and in this case we have that the parameterization $P$ is formally well defined to all orders. 
That is, Equation \eqref{eq:invEq} has a well defined formal series solution satisfying the first order 
constraints of Equations \eqref{eq:constraint1} and \eqref{eq:constraint2}.  
\end{remark}  

\begin{remark}[Uniqueness up to rescaling of the first order data]
The unique solvability of the homological equations,
assuming non-resonance of the unstable eigenvalues, 
gives that the solution $P$ at $u_0$ is unique up to the choice of the 
scalings of the eigenfunctions.  The choice of the scaling of the eigenfunctions 
directly effects the decay of the coefficients $p_{n_1, \ldots, n_M}$ as discussed 
in \cite{MR3437754,parmPDEone}.  For this reason we always fix the 
domain of the parameterization to be $\mathbb{B} = [-1,1]^M$, and choose the 
scaling of the eigenvectors so that the coefficients decay rapidly.
Of course while choosing smaller scalings for the eigenvectors 
provides faster coefficient decay, it also means that the image of $\mathbb{B}$
is smaller in $\mathcal{H}$.  That is, smaller scalings stabilize the numerics 
but reveal a smaller portion of the local unstable manifold.  In practice we 
must strike a balance between the polynomial order of the calculation 
(at what order do we truncate the formal series?) the scaling of the eigenvectors 
and the size of the local unstable manifold we compute.   
\end{remark}

\subsection{Automatic differentiation of power series} \label{sec:autoDiff}
A critical step in any explicit example is to work out the dependence of the coefficients
$q_{n_1, \ldots, n_M}$ of the nonlinear composition on the unknown 
coefficients $p_{n_1, \ldots, n_M}$.
This is essential for defining the right hand side $S_{n_1, \ldots, n_M}$ of the 
Homological equation \eqref{eq:homEqHighLevel}.
This challenge reduces to repeated application of the Cauchy product formula whenever
$f(\cdot, x,y)$ has polynomial nonlinearity.  

For example consider the case where $f$ is a quadratic function of the form
\[
f(u, x,y) = a(x,y) u^2.
\]  
Then 
\[
f(P(\theta, x,y), x,y)
\]
\begin{align*}
 &= a(x,y) \left(
\sum_{n_1 = 0}^\infty \ldots \sum_{n_M = 0}^\infty p_{n_1, \ldots, n_M}(x,y) \theta_1^{n_1} \ldots \theta_M^{n_M}
\right) \left( 
\sum_{n_1 = 0}^\infty \ldots \sum_{n_M = 0}^\infty p_{n_1, \ldots, n_M}(x,y) \theta_1^{n_1} \ldots \theta_M^{n_M}
\right)  \\
&= a(x,y) \sum_{n_1 = 0}^\infty \ldots \sum_{n_M = 0}^\infty \left(
\sum_{k_1 = 0}^{n_1} \ldots \sum_{k_M = 0}^{n_M} 
p_{n_1 - k_1, \ldots, n_M - k_M}(x,u) p_{k_1, \ldots, k_M}(x,y) \right)
\theta_1^{n_1} \ldots \theta_M^{n_M} \\
& =  \sum_{n_1 = 0}^\infty \ldots \sum_{n_M = 0}^\infty \left(
\sum_{k_1 = 0}^{n_1} \ldots \sum_{k_M = 0}^{n_M} 
a(x,y) p_{n_1 - k_1, \ldots, n_M - k_M}(x,y) p_{k_1, \ldots, k_M}(x,y) \right)
\theta_1^{n_1} \ldots \theta_M^{n_M} \\
& = \sum_{k_1 = 0}^{n_1} \ldots \sum_{k_M = 0}^{n_M}  q_{n_1, \ldots, n_M}(x,y) \theta_1^{n_1} \ldots \theta_M^{n_M},
\end{align*}
and we see that 
\begin{align*}
 q_{n_1, \ldots, n_M}(x,y) &= \sum_{k_1 = 0}^{n_1} \ldots \sum_{k_M = 0}^{n_M} 
a(x,y) p_{n_1 - k_1, \ldots, n_M - k_M}(x,y) p_{k_1, \ldots, k_M}(x,y)  \\
&= 2 a(x,y) p_{0, \ldots, 0}(x,y) p_{n_1, \ldots, n_M}(x,y) + \mbox{``lower order terms''} \\
&= 2 \frac{\partial}{\partial u} f(u_0, x,y) p_{n_1, \ldots, n_M}(x,y) + \mbox{``lower order terms''},
\end{align*}
as promised above.
Indeed the ``lower order terms'' have the explicit form 
\[
S_{n_1, \ldots, n_M} = \sum_{k_1 = 0}^{n_1} \ldots \sum_{k_M = 0}^{n_M} 
\hat{\delta}_{n_1, \ldots, n_M}^{k_1, \ldots, k_M} 
a(x,y) p_{n_1 - k_1, \ldots, n_M - k_M}(x,y) p_{k_1, \ldots, k_M}(x,y) 
\]
where the coefficient 
\[
\hat{\delta}_{n_1, \ldots, n_M}^{k_1, \ldots, k_M} = 
\begin{cases}
0 & \mbox{if } k_1 = \ldots = k_M = 0 \\
0 & \mbox{if } k_1 = n_1, \ldots, k_M = n_M \\
1 & \mbox{otherwise}
\end{cases},
\]
appears in the sum to indicate that both of the terms with
$p_{n_1, \ldots, n_M}(x,y)$ have been removed.  

When $f$ contains non-polynomial terms, calculating the 
$q_{n_1, \ldots, n_M}$ is more delicate.  We employ a semi-numerical 
technique based on the idea that many typical nonlinearities  
appearing in applications are themselves solutions of polynomial 
differential equations.  This is exploited in fast recursion 
schemes.

Consider for example the case of 
\[
f(u, x,y) = a(x,y) e^{-u}.
\]
Let 
\begin{equation*}  
 P(\theta, x, y) = \sum_{n_1 = 0}^\infty \ldots \sum_{n_M = 0}^\infty p_{n_1, \ldots, n_M}(x,y) \theta_1^{n_1} \ldots \theta_M^{n_M},
\end{equation*}
and write  
\begin{equation} \label{eq:fP_inRadGrad}
Q(\theta, x,y) = 
\sum_{n_1 = 0}^\infty \ldots \sum_{n_M = 0}^\infty q_{n_1, \ldots, n_M}(x,y) \theta_1^{n_1} \ldots \theta_M^{n_M} = f(P(\theta, x, y)).
\end{equation}
The following idea is described in detail in Chapter $2$ of \cite{MR3467671}.  We apply the 
\textit{radial gradient} -- the first order partial differential operator given by 
\[
\nabla_{\theta} = \theta_1 \frac{\partial}{\partial \theta_1} + 
\ldots + 
\theta_M \frac{\partial}{\partial \theta_M}, 
\]
to both sides of Equation \eqref{eq:fP_inRadGrad} and obtain 
\[
\nabla_\theta f(P(\theta, x,y), x,y) = \nabla_\theta Q(\theta, x,y).
\]
That is
\[
\nabla_\theta f(P(\theta, x,y), x,y) 
\]
\begin{align*}
&= \theta_1\frac{\partial}{\partial u} f(u,x,y)\left|_{u = P(\theta, x,y)} \right. \frac{\partial}{\partial \theta_1} P(\theta, x,y) 
+ \ldots + 
\theta_M\frac{\partial}{\partial u} f(u,x,y)\left|_{u = P(\theta, x,y)}  \right. \frac{\partial}{\partial \theta_M} P(\theta, x,y) \\
&= -a(x,y) e^{- P(\theta, x,y)} \left(
\theta_1 \frac{\partial}{\partial \theta_1} P(\theta, x,y) 
+ \ldots + 
\theta_M \frac{\partial}{\partial \theta_M} P(\theta, x,y) 
\right) \\
& = - Q(\theta, x,y) \nabla_\theta P(\theta, x,y)
\end{align*}
\[
 = - \left(
 \sum_{n_1 = 0}^\infty \ldots \sum_{n_M = 0}^\infty 
q_{n_1, \ldots, n_M}(x,y) \theta_1^{n_1} \ldots \theta_M^{n_M}
\right)
\]
\[
 \left(
 \sum_{n_1 = 0}^\infty \ldots \sum_{n_M = 0}^\infty 
(n_1 + \ldots + n_M) p_{n_1, \ldots, n_M}(x,y) \theta_1^{n_1} \ldots \theta_M^{n_M}
\right)
\]
\[
= -  \sum_{n_1 = 0}^\infty \ldots \sum_{n_M = 0}^\infty \left(
\sum_{k_1 =0}^{n_1} \ldots \sum_{k_M =0}^{n_M} 
(k_1 + \ldots + k_M)q_{n_1 - k_1, \ldots, n_M - k_M} p_{k_1, \ldots, k_M} \right)
 \theta_1^{n_1} \ldots \theta_M^{n_M},
\]
on the left, and 
\[
\nabla_\theta Q(\theta, x,y) =  \sum_{n_1 = 0}^\infty \ldots \sum_{n_M = 0}^\infty 
(n_1 + \ldots + n_M) q_{n_1, \ldots, n_M}(x,y) \theta_1^{n_1} \ldots \theta_M^{n_M}
\]
on the right. Matching like powers and isolating $q_{n_1, \ldots, n_M}$ leads to 
\[
q_{n_1, \ldots, n_M} = \frac{-1}{n_1 + \ldots + n_M} 
\sum_{k_1 =0}^{n_1} \ldots \sum_{k_M =0}^{n_M} 
(k_1 + \ldots + k_M)q_{n_1 - k_1, \ldots, n_M - k_M} p_{k_1, \ldots, k_M}.
\] 
Then the complexity of computing the power series coefficients of $a(x,y) e^{-P(\theta, x,y)}$ 
is the complexity of a single Cauchy product.  The additional cost is that the coefficients 
of $Q$ have to be stored in addition to those of $P$.

Such methods for formal series manipulations are referred to by many authors as
\textit{automatic differentiation for power series}, and they facilitate rapid computation 
of the formal series coefficients of compositions with all the elementary functions.  
  A classic reference which includes 
an in depth historical discussion is found in Chapter $4$, Section $6$ of \cite{MR633878}.
See also the discussion of software implementations found in \cite{MR2146523}.

\section{Applications} \label{sec:applications}

\subsection{A first worked example: Fisher's Equation } \label{sec:composedHomEq}

Consider the parabolic PDE 
\[
\frac{\partial}{\partial t}u=\Delta u +\alpha u (1-u),
\]
on the $\mathbb{L}$ domain $\Omega$ illustrated in the left-most frame
of Figure \ref{fig:planarDomains}, subject to the 
Neumann boundary conditions 
\[
\nabla u \cdot \textbf{n}|_{\partial \Omega}=0.
\]
Here \textbf{n} is a unit vector normal to $\partial\Omega$.
This reaction-diffusion equation was introduced by Ronald Fisher in the context of population 
dynamics, as a toy model for the propagation of advantageous genes. 
Letting 
\[
F(u) = \Delta u +\alpha u (1-u),
\]
we see that the problem describes an evolution equation as
in Equation \eqref{eq:evolution}.

Recall that an equilibrium solution has $F(u) = 0$, and note that 
$0$ is always an equilibrium.  We refer to $0$ as the homogeneous 
background solution, and note that   
while for small $\alpha$ it is stable, it 
looses stability as $\alpha$ increases.
Each time an eigenvalue of the homogeneous solution crosses the 
imaginary axis, the bifurcation gives rise to a pair of non-trivial equilibrium
solutions.    The first pair of non-trivial equilibria to appear are stable initially, 
but loose stability as $\alpha$ is further increased.  Hence, at 
$\alpha=2.7$ we can find a non-trivial equilibrium solution with Morse index 1, and 
Morse index 2 when $\alpha = 9$.  These equilibrium solutions 
have one and two dimensional attached unstable
manifolds.   In the remainder of this section we discuss in detail the parameterization of the 
two dimensional unstable manifold for this otherwise simple example.

To find equilibria, we study the nonlinear elliptic BVP
\[
F(u) = \Delta u + \alpha u (1-u) = 0,
\]
subject to the same natural boundary conditions on $\Omega$.
The weak formulation is
\[
\mathcal{F}(u)\phi=-\int_{\Omega} \nabla u \cdot \nabla \phi +\int_{\Omega} \alpha u(1-u)\phi =0,
\]
and, using the notation of Section \ref{sec:finiteElements},
triangulate $\Omega$ and solve for the coefficients of the finite 
element representation 
$u_h=\sum\limits_{j=1}^{{nb}} c_j \phi_j$ of $u$.
In order to construct this projection, 
define the linear basis functions $\phi_j$ as
\[ 
\phi_j(n_i)=
\begin{cases} 
1 & j = i \\
0 & j \ne i 
\end{cases},
\]
where $n_i$ denotes the $i-th$ vertex in the triangulation. 
Note that in this case, ${nb}={nn}$.
Letting $\phi=\phi_i$ for $1\leq i\leq {nb}$ leads to the nonlinear system 
of ${nb}$ equations in ${nb}$ unknowns, given by 
\[
\mathcal{F}^{h}_i(c)=-\int_{\Omega} \left(\sum_{j=1}^{{nb}} c_j \nabla \phi_j \right) \cdot \nabla \phi_i +  \int_{\Omega} \alpha \left(\sum_{j=1}^{{nb}} c_j \phi_j \right) \left(1-\sum_{j=1}^{{nb}} c_j \phi_j \right)\phi_i =0,
\]
which we solve using the Newton's Method (for $c=(c_1,c_2,...,c_{nb})$). 
More precisely, let 
$\mathcal{F}^{h}(c)=(\mathcal{F}^{h}_1(c),...,\mathcal{F}^{h}_{{nn}}(c))=(\mathcal{F}(u_h)\phi_1,...,\mathcal{F}(u_h)\phi_{nn})$.  The $k$'th Newton's step is given by
\[
c^{(k)}= c^{(k-1)} -{D\mathcal{F}^{h}\left(c^{(k-1)}\right)}^{-1}\mathcal{F}^{h}\left(c^{(k-1)}\right),
\]
where $u_h^{(k)}= \sum_{j=1}^{{nb}} c^{(k)}_j \phi_j$ and $D\mathcal{F}^{h}(c)= -\left(\int_{\Omega} \nabla \phi_j \cdot \nabla \phi_i\right) + \left(\int_{\Omega} \frac{ \partial N(c)}{ \partial c_j} \phi_i\right)$.

\begin{figure}[!t] 
\begin{center}
	\begin{minipage}[b]{0.32\textwidth}
		\includegraphics[width=\textwidth]{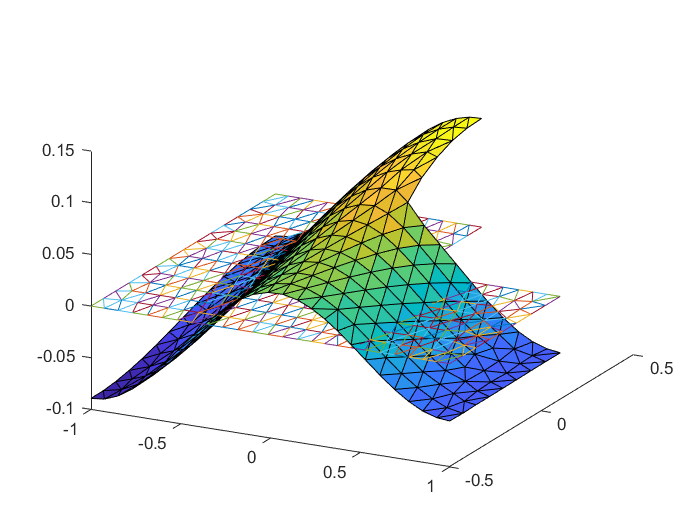}
	\end{minipage}
    \begin{minipage}[b]{0.32\textwidth}
	\includegraphics[width=\textwidth]{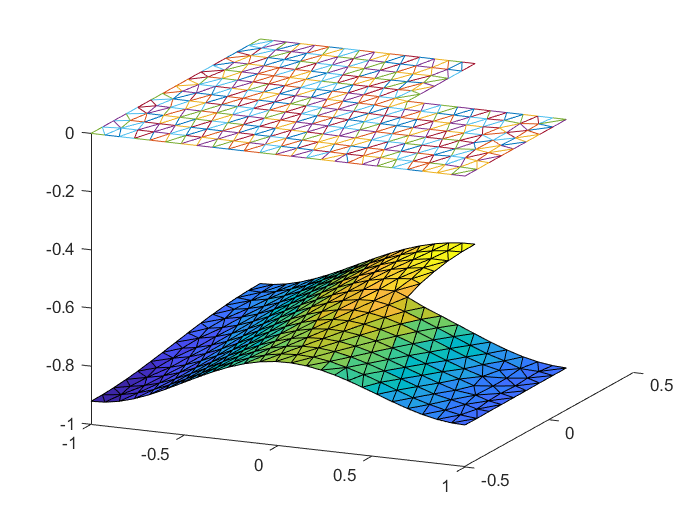}
	\end{minipage}
	\begin{minipage}[b]{0.32\textwidth}
		\includegraphics[width=\textwidth]{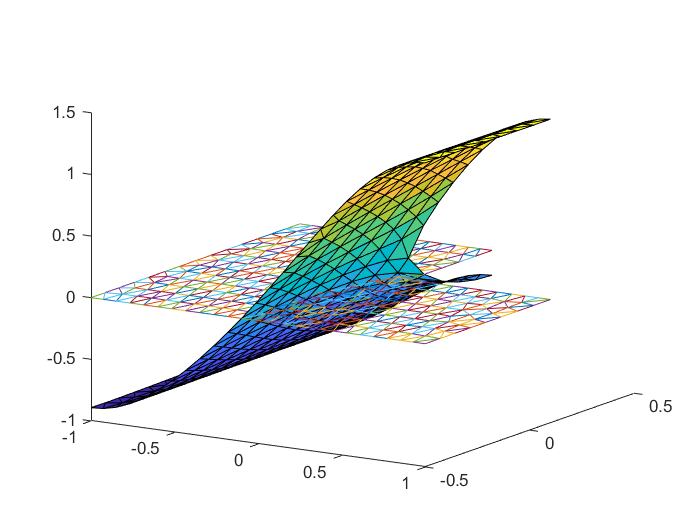}
	\end{minipage}
    \end{center}
    \caption{ Fisher's equation with $\alpha=9$, $ne=515$. 
    Left: Equilibrium solution, with Morse index 2.
    Center: Eigenfunction for $\lambda_1=9.04$.
    Right:  Eigenfunction for $\lambda_2=7.16$.}
\end{figure}


Once the approximate solution $u_0$ is computed we proceed to solve the eigenvalue-eigenvector problem
\[
\Delta \xi + \alpha (1-2 u_0) \xi -\lambda \xi =0.
\]
That is, we compute the projection $\xi_h=\sum\limits_{j=1}^{{nb}} c_j \phi_j$ 
in the weak formulation, which leads to 
\[
-\int_{\Omega} \left(\sum_{j=1}^{{nb}} c_j \nabla \phi_j \right) \cdot \nabla \phi_i +  \int_{\Omega} \alpha (1-2u_0) \left(\sum_{j=1}^{{nb}} c_j \phi_j \right)\phi_i =\int_{\Omega} \lambda \left(\sum_{j=1}^{{nb}} c_j \phi_j \right) \phi_i  
\]
or
\[
\Big(- \int_{\Omega} \nabla \phi_j \cdot \nabla \phi_i + \alpha (1-2u_0)\phi_j \phi_i\Big)c = \lambda \Big(\int_{\Omega}\phi_j \phi_i\Big)c.
\]


After computing the unstable eigenvalues $\lambda_1$ and $\lambda_2$ and the 
associated eigenfunctions
$\xi_1$ and $\xi_2$, we proceed to solve the invariance equation \eqref{eq:invEq} 
specialized to the present situation.  That is, we consider the weak form of the equation  
\[
F(P(\theta))=\lambda_1 \theta_1 \frac{\partial }{\partial \theta_1} P(\theta)
+\lambda_2 \theta_2 \frac{\partial }{\partial \theta_2} P(\theta),
\]
where
 \[
 P(\theta)=\sum\limits_{m=0}^{\infty}\sum\limits_{n=0}^{\infty}p_{m,n}(x,y)\theta_1^m \theta_2^n,
\] 
with $p_{0,0}=u_0$, $p_{1,0}=\xi_1$ and $p_{0,1}=\xi_2$. 
Taking the projection $p_{m,n}=\sum\limits_{j=1}^{{nb}}c_j^{(m,n)} \phi_j$, leads to 
\[
\Big(- \int_{\Omega} \nabla \phi_j \cdot \nabla \phi_i + \alpha( 1-2u_0- m\lambda_1 -n\lambda_2 ) \phi_j \phi_i\Big) \Big( c_i^{(m,n)}\Big) =  \Big(\int_{\Omega} s_{(m,n)} \phi_i \Big),
\]
for $m +n\geq 2$, 
which is 
\[
\Big(D\mathcal{F}^{h}(c^{(0)})-(\lambda_1 m + \lambda_2 n) \int_{\Omega} \phi_j \phi_i \Big) c^{(m,n)}=\Big( \int_{\Omega}  s_{(m,n)} \phi_i \Big),
\]  
with
\[
s_{(m,n)}=\alpha \sum\limits_{i=0}^{m} \sum\limits_{j=0}^{n} \delta(i,j) p_{i,j}p_{m-i,n-j},
\]
and
\[ 
\delta(i,j)=
\begin{cases} 
0 & (i,j)=(0,0) \hspace{2mm} \text{or} \hspace{2mm} (i,j)=(m,n) \\
1 & \text{otherwise}
\end{cases}.
\]

\begin{figure}[!t]
\begin{center}
	\begin{minipage}[l]{0.49\textwidth}
		\includegraphics[width=\textwidth]{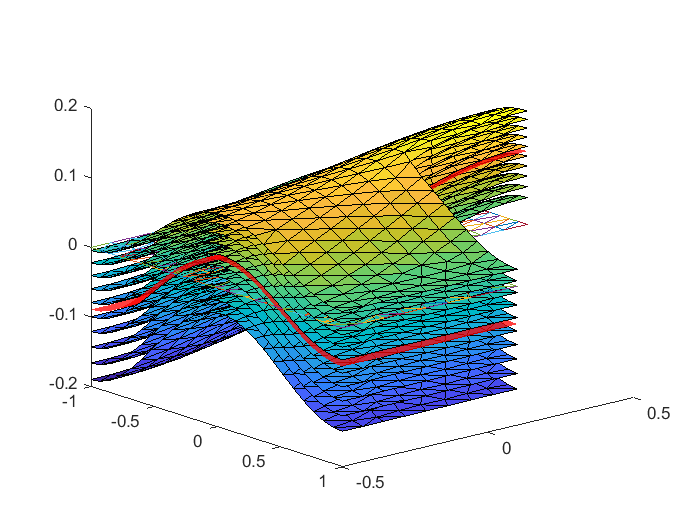}
	\end{minipage}
    \begin{minipage}[r]{0.49\textwidth}
	\includegraphics[width=\textwidth]{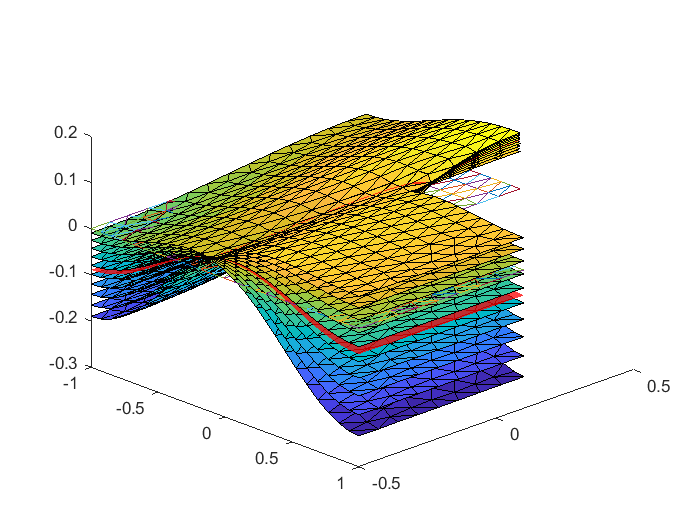}
	\end{minipage}
    \end{center}
    \caption{Left: 10 functions on the fast manifold approximated to order $N=30$ with Invariance equation error of 1.34e-10 with respect to the $L^2$ norm. 
    Right: 10 functions on the slow manifold approximated to order $N=30$ with Invariance equation error of 4.66e-08 with respect to the $L^2$ norm.}
\label{fig:fisher2dMani}
\end{figure}

As anticipated in Section 3.1, the homological equations are linear 
elliptic PDEs, and we solve them 
recursively to any desired order using the Finite Element Method. 
Figure \ref{fig:fisher2dMani} shows a few functions in the 
\textit{fast manifold} (1d manifold associated to the largest positive eigenvalue) 
and \textit{slow manifold} (1d manifold associated to the smallest positive 
eigenfunction) approximated up to order $N=30$. 
 
\begin{figure}[!t] 
\begin{center}
	\begin{minipage}[b]{0.32\textwidth}
		\includegraphics[width=\textwidth]{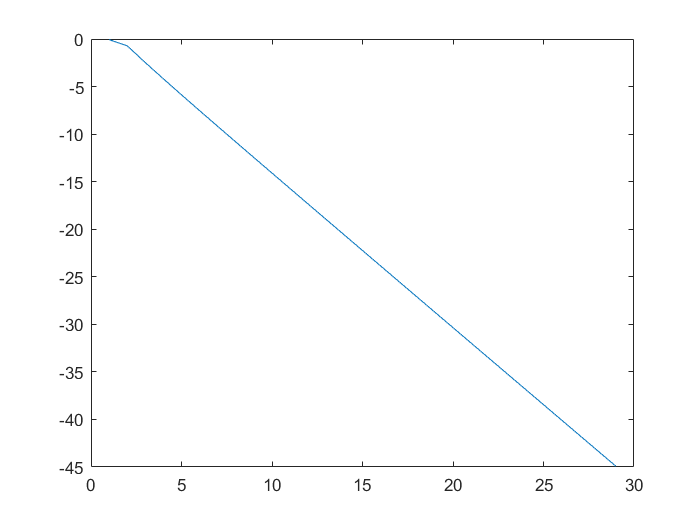}
	\end{minipage}
    \begin{minipage}[b]{0.32\textwidth}
	\includegraphics[width=\textwidth]{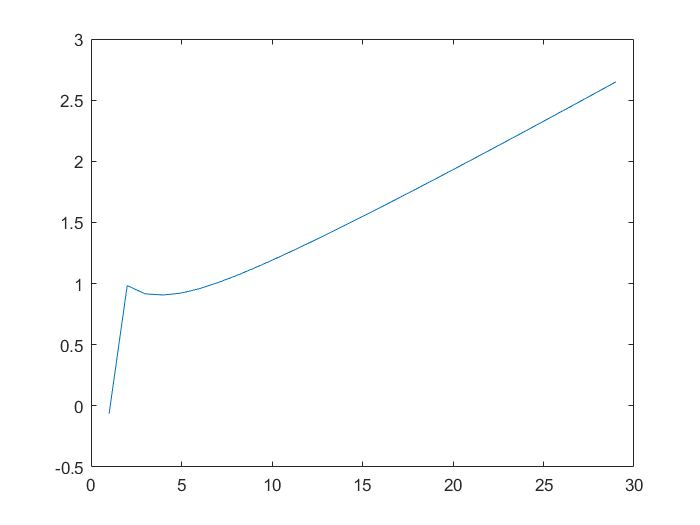}
	\end{minipage}
	\begin{minipage}[b]{0.32\textwidth}
		\includegraphics[width=\textwidth]{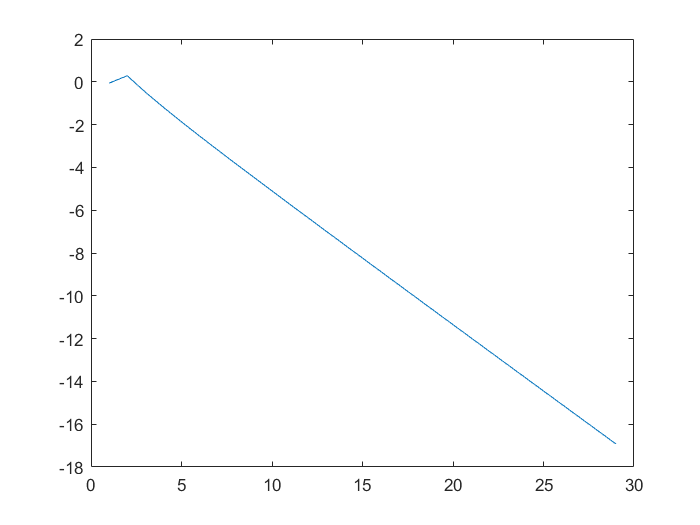}
	\end{minipage}
    \end{center}
    \caption{\textbf{Coefficient growth:} three plots of the magnitude of the parameterization coefficients
    as a function of the order of the coefficients. (Horizontal axis is the order of the coefficient and vertical 
    axis is the base ten logarithm $L^2$ norm of the coefficient function).
    Left: The scaling of the eigenvector is too small, and the coefficients decay too fast.  Coefficients after order then
     are below machine precision in $L^2$ norm (smaller than $10^{-16}$) and hence do not  contribute 
     significantly to the accuracy of the polynomial approximation.  
    Center: The eigenvector scaling is choosen too large, and 
    now the $p_{m,n}$'s grow exponentially fast.  This introduces numerical instabilities into the approximation.
    Right: The scaling is chosen just right: they decay exponentially fast at a rate chosen so that the 
    $N$-th order coefficients reach machine precision.}
\label{fig:coeffDecay}
\end{figure}

The effect of the scaling of the eigenvectors on the decay of the coefficients is illustrated
 in Figure \ref{fig:coeffDecay} .

\subsection{A reaction diffusion equation with non-polynomial nonlinearity: one unstable eigenvalues}

In this section we derive the homological equations for a non-polynomial problem.
We consider the reaction diffusion equation with Ricker type exponential nonlinearity 
given by 
\begin{equation} \label{eq:FisherRicker}
u_t=\Delta u + \alpha u \left(0.5-e^{-u} \right). 
\end{equation}
We refer to this problem as the Fisher-Ricker (FR) equation,
and take parameter $\alpha = -4.7$. 
Letting 
\[
F(u) = \Delta u + \alpha u \left(0.5-e^{-u} \right),
\]
we obtain an evolution equation of the kind given in 
Equation \eqref{eq:evolution}.

To find the equilibrium solution consider the weak form of the 
equation $F(u) = 0$, project into a finite element space of
piecewise linear functions, and solve
\[
\mathcal{F}^{h}_i(c)=
-\int_{\Omega} \left(\sum_{j=1}^{{nb}} c_j \nabla \phi_j \right) \cdot \nabla \phi_i +  \int_{\Omega} \alpha 
\left(\sum_{j=1}^{{nb}} \phi_j\right) \left(0.5- \exp\{-\sum_{j=1}^{{nb}} c_j \phi_j\}\right) \phi_i  =0.
\]
The corresponding eigenvalue-eigenfunction problem is 
\[
D\mathcal{F}^{h}(c^{(0)})c = \lambda\Big( \int_{\Omega} \phi_j \phi_i \Big)c.
\] 

Suppose now that $u_0$ is an equilibrium solution with Morse index 1, let
$\lambda$ denote the unstable eigenvalue, and $\xi$ be a corresponding eigenfunction. 
We seek a parameterization of the form
$P(\theta)=\sum\limits_{n=0}^{\infty}p_n \theta^n$ 
solving the 1D Invariance Equation
\[
F(P(\theta)) = \theta \lambda \frac{d}{d \theta} P(\theta), 
\]
which, after expanding $P(\theta)$ as a power series becomes 
\[
\sum\limits_{n=0}^{\infty}\Delta p_n \theta^n +\alpha \left(
\sum\limits_{n=0}^{\infty}p_n \theta^n\right) \left(0.5-
\exp \left(-\sum\limits_{n=0}^{\infty}p_n \theta^n \right)
\right)=\lambda\sum\limits_{n=0}^{\infty} n p_n \theta^n.
\]
Here, the $p_n = p_n(x,y)$ are functions defined on $\Omega$
satisfying the boundary conditions.

The challenge is to compute the power series expansion of  the exponential. 
To this end, we introduce the new variable 
\[
Q(\theta):= e^{-P(\theta)}=\sum\limits_{n=0}^{\infty} q_n \theta^n,
\]
and apply the automatic differentiation technique described in Section \ref{sec:autoDiff}. 
That is, we note that 
$Q'=-QP'$, 
and expand the relation as a product of power series.
Matching like powers, we obtain
\[
(n+1)q_{n+1}=-\sum\limits_{j=0}^{n} (j+1) p_{j+1} q_{n-j},
\]
and isolating the $n$-th order terms we have
\begin{equation} \label{eq:for_qn}
q_n=-p_n q_0 - \frac{1}{n} \sum\limits_{j=0}^{n-2}(j+1)p_{j+1}q_{n-1-j}.
\end{equation}
Note that Equation \eqref{eq:for_qn} involves only sums and products
of the functions $p_i(x,y), q_j(x,y)$, for $0 \leq i,j \leq n$, and that these 
operations are well defined for $p_n, q_n$ in any finite element space.
Equation \eqref{eq:for_qn} then 
allows us to compute $q_n$ to any desired order, assuming that 
$p_n, \ldots, p_0$, and $q_{n-1}, \ldots, q_0$ are known.  

Returning to the Invariance Equation and using the recursive formula for $q_n$ we obtain
that for $n \geq 2$, the $p_n$ solve 
\[
\Delta p_n +\alpha( 0.5-q_0-\lambda n )p_n -\alpha p_0 q_n = \alpha \sum\limits_{j=1}^{n-1} p_j q_{n-j},
\]
or 
\[
\Delta p_n +\alpha(0.5-q_0+p_0q_0- \lambda n)p_n= s_n,
\]
where
\[
 s_n=\alpha \sum\limits_{j=1}^{n-1} p_j q_{n-j} - \frac{\alpha p_0}{n}
\sum\limits_{j=0}^{n-2} (j+1) p_{j+1}q_{n-1-j}.
\]
Passing to the weak form, we find that the 
coefficients $p_{n}=\sum\limits_{j=1}^{{nb}}c_j^{(n)} \phi_j$ solve
the homological equations
\begin{equation} \label{eq:expHom}
\Big(D\mathcal{F}^{{h}}(c^{(0)})-\lambda n \int_{\Omega} \phi_j \phi_i \Big) c^{(n)}
=\Big( \int_{\Omega}  s_n \phi_i \Big),
\end{equation}
for $n \geq 2$.
Notice that $s_n$ only depends on $p_k$'s and $q_k$'s with $0 < k<n$. 
Then if $p_0, \ldots, p_{n-1}$ and $q_0, \ldots, q_{n-1}$
are known, $p_n$ is computed by solving Equation \eqref{eq:expHom}.
Once $p_n$ is known, we update Equation \eqref{eq:for_qn}
to obtain $q_n$.

\subsection{A reaction diffusion equation with non-polynomial nonlinearity: 
two unstable eigenvalues}

A modification of the method just discussed allows us to compute 
higher dimensional manifolds in problems with non-polynomial nonlinearities.
Consider again Equation \eqref{eq:FisherRicker}, 

this time with $\alpha=-4.41$.  At this parameter value there is a non-trivial 
equilibrium $u_0$ with Morse index 2.
Let $\lambda_1$ and $\lambda_2$ denote the unstable eigenvalues
and $\xi_1$, $\xi_2$ denote an associated pair of unstable eigenfunctions.

\begin{figure}[!t]
\begin{center}
	\begin{minipage}[l]{0.32\textwidth}
		\includegraphics[width=\textwidth]{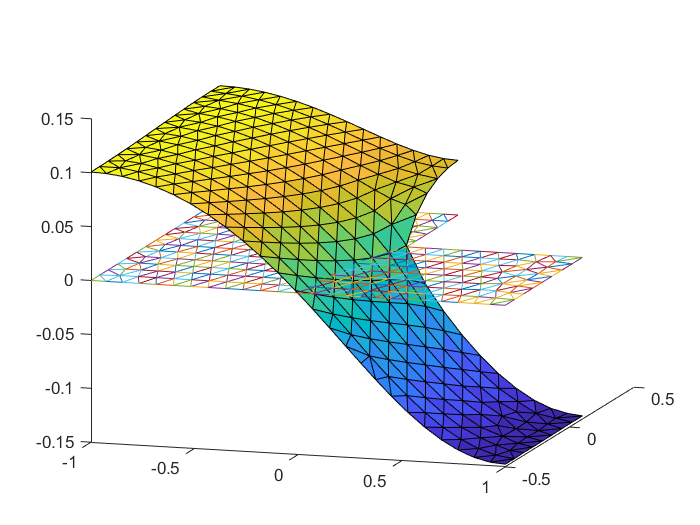}
	\end{minipage}
    \begin{minipage}[r]{0.32\textwidth}
	\includegraphics[width=\textwidth]{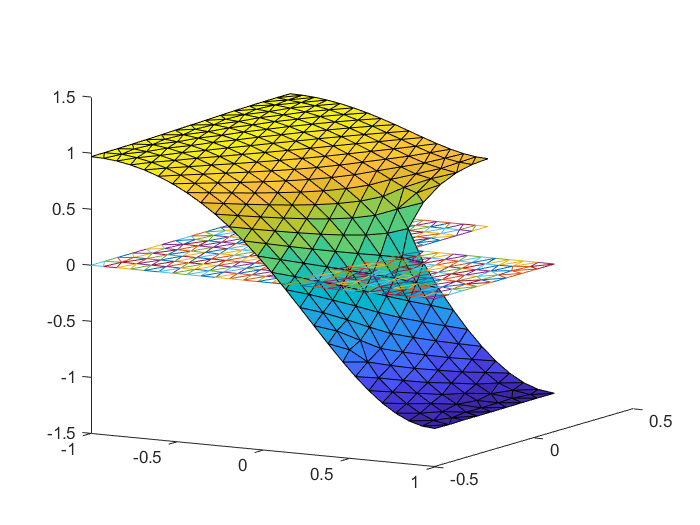}
	\end{minipage}
	\begin{minipage}[l]{0.32\textwidth}
		\includegraphics[width=\textwidth]{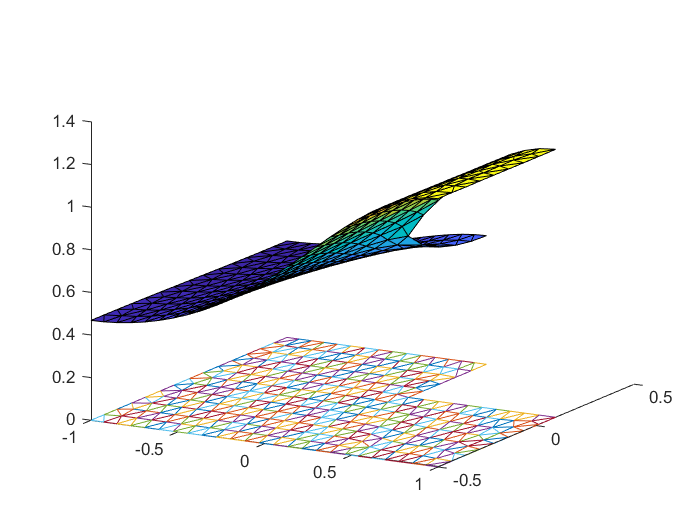}
	\end{minipage}
    \end{center}
    \caption{Fisher-Ricker equation with $\alpha=-4.7$, $ne=515$. Left: Equilibrium solution. 
    Center:  Eigenfunction $\xi_1$ with $\lambda_1=2.41$.
    Right: Eigenfunction $\xi_2$  with $\lambda_2=0.05$.}
\end{figure}

Recall that for an equilibrium with Morse index 2, the invariance equation becomes 
\[
F(P(\theta))=\lambda_1 \theta_1 \frac{\partial }{\partial \theta_1} P(\theta)
+\lambda_2 \theta_2 \frac{\partial }{\partial \theta_2} P(\theta),
\]
and we seek a power series solution of the form
\[
P(\theta)=\sum\limits_{m=0}^{\infty}\sum\limits_{n=0}^{\infty}p_{m,n}(x,y)\theta_1^m \theta_2^n,
\]
with 
\[
p_{00} = u_0, \quad \quad \quad p_{10} = \xi_1, \quad \mbox{and} \quad 
p_{01} = \xi_2, 
\]
and where $p_{m,n}$ for $m+n \geq 2$ are to be determined. 
To work out the exponential nonlinearity, define the auxiliary equation 
\[
Q:=\exp \left(-P(\theta)\right)=
\sum\limits_{m=0}^{\infty}\sum\limits_{n=0}^{\infty}q_{m,n}(x,y)\theta_1^m \theta_2^n.
\]
Taking the radial gradient of both sides of this equation,
as discussed in Section \ref{sec:autoDiff}, leads to 
\[
\nabla_\theta Q(\theta) = \nabla_\theta \left( \exp^{-P(\theta)} \right),
\]
or 
\[
\theta_1 \frac{\partial }{\partial \theta_1} Q
+\theta_2 \frac{\partial }{\partial \theta_2} Q = -Q\left( \theta_1 \frac{\partial }{\partial \theta_1} P
+ \theta_2 \frac{\partial }{\partial \theta_2} P \right).
\]
Plugging in the power series, computing the derivatives (formally), and 
matching like powers leads to 
\[
\sum\limits_{m,n\geq1}(m+n)q_{m,n}\theta_1^m\theta_2^n = -\Big(\sum\limits_{m,n\geq1}(m+n)p_{m,n}\theta_1^m\theta_2^n\Big)\Big(\sum\limits_{m,n\geq0}q_{m,n}\theta_1^m\theta_2^n\Big).
\]
Expanding the Cauchy products, and isolating $q_{m,n}$ leads to 
\[
q_{m,n}=-\frac{1}{(m+n)} \sum\limits_{i=1}^{m}\sum\limits_{j=1}^{n}(i+j)p_{i,j}q_{m-i,n-j}.
\]
Note that this requires only additions and multiplications, all well defined operations 
for finite element basis functions.  

\begin{figure}[!t] \label{fig:expMorse2}
\begin{center}
    \begin{minipage}[r]{0.49\textwidth}
	\includegraphics[width=\textwidth]{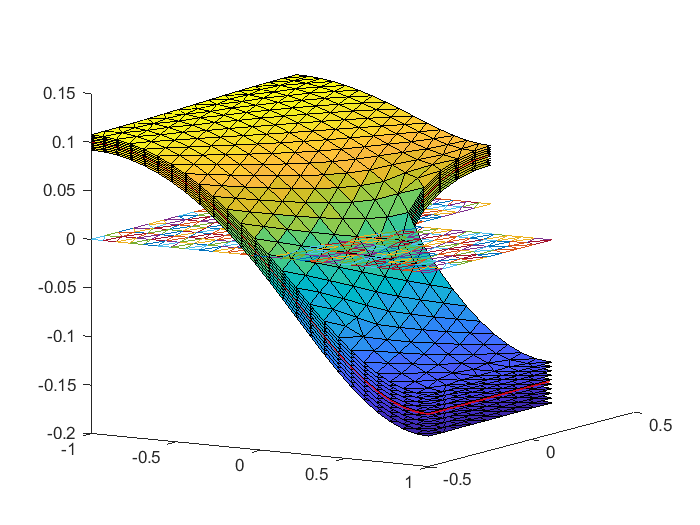}
	\end{minipage}
	\begin{minipage}[l]{0.49\textwidth}
		\includegraphics[width=\textwidth]{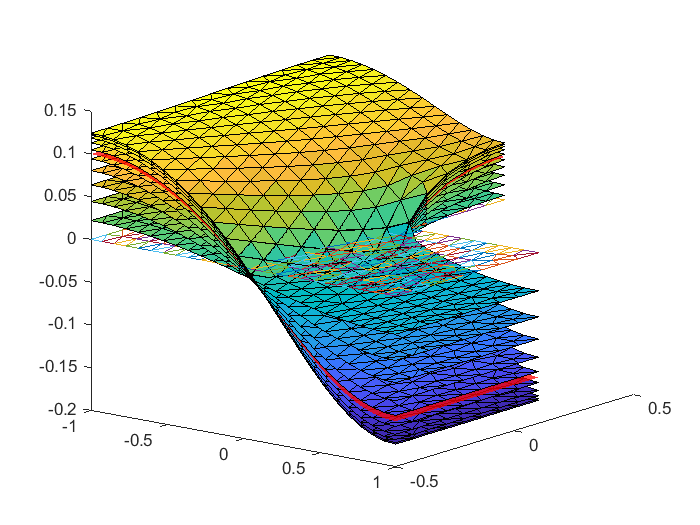}
	\end{minipage}
    \end{center}
    \caption{    	
    Left: 10 functions on the fast manifold approximated to order $N=30$ with Invariance equation error of 9.21e-10 with respect to the $L^2$ norm. 
    Right: 10 functions on the slow manifold approximated to order $N=30$ with Invariance equation error of 6.12e-08 with respect to the $L^2$ norm.}
    
\end{figure}

Returning to the Invariance Equation and using the recursive formula for $q_{m,n}$ we have 
\[
\Delta p_{m,n} +\alpha( 0.5-q_{0,0}-\lambda_1 m -\lambda_2 n )p_{m,n} -\alpha p_{0,0} q_{m,n}= \alpha \sum\limits_{i=0}^{m}\sum\limits_{j=0}^{n}q_{i,j}p_{m-i,n-j} \delta(i,j),
\]
so that the strong form of the homological equation is 
\[
\Delta p_{m,n} +\alpha(0.5-q_{0,0} + p_{0,0}q_{0,0} - \lambda_1 m - \lambda_2 n)p_{m,n}= s_{m,n},
\]
with 
\[
s_{m,n}=\alpha \sum\limits_{i=0}^{m}\sum\limits_{j=0}^{n}q_{i,j}p_{m-i,n-j} \delta(i,j) 
- \frac{\alpha p_{0,0}}{(m+n)}\sum\limits_{i=1}^{m-1}\sum\limits_{j=1}^{n-1} (i+j) p_{i,j}q_{m-i,n-j},
\]
a linear, elliptic BVP for each $m + n \geq 2$ as desired. 
Passing to the weak form leads to  
\[
\Big(D\mathcal{F}^{h}(c^{(0)})-(\lambda_1 m +\lambda_2 n) \int_{\Omega} \phi_j \phi_i \Big) c^{(m,n)}
=\Big( \int_{\Omega}  s_n \phi_i \Big),
\] 
which we solve recursively via the finite element method, obtaining the 
parameterization coefficients to any desired order (updating the equation for 
$q_{mn}$ as we go).  Results are illustrated in Figure \ref{fig:expMorse2}.

\subsection{Higher order PDEs: a Kuramoto-Sivashinsky small term} \label{sec:KS_manifolds}

We now consider a higher order problem, whose leading diffusion term is given by 
the biharmonic Laplacian.  The biharmonic operator often 
appears in models of thin structures that react elastically to external forces.
Consider the Kuramoto-Sivashinsky equation given by 
\[
F(u)=-\Delta^2 u - \Delta u -0.5|\nabla u|^2,
\]
\[
u|_{\partial \Omega}=0  \quad \quad \quad  \nabla u \cdot \textbf{n}|_{\partial \Omega}=0.
\]
which models the propagation of a flame front and it is known to exhibit chaotic dynamics.
We refer to \cite{MR2496834,MR1815526,MR3390404,MR796268,KT}
for more complete discussion of the equation, it's physics, and it's dynamical properties. 
 
Since the differential operator is fourth order, higher order Finite Elements are required. 
The purpose of this section is to illustrate the use of the 
parameterization method in a higher order problem. To exhibit our approach in a simple manner,
we start with an already computed solution of Fisher and introduce a biharmonic term and nonlinearity 
as a perturbation with natural boundary conditions.

Specifically, we take 
\[
F_0(u)= \Delta u + \alpha u(1-u)
\quad \quad \quad \nabla u \cdot \textbf{n}|_{\partial \Omega}=0,
\]
with weak formulation
\[
\mathcal{F}_{0}(u) \phi = -\int_{\Omega} \nabla u \cdot \nabla \phi +\int_{\Omega} \alpha u(1-u)\phi =0, 
\]

and let $\mathcal{F}_{\epsilon}(u) \in {H^2}^{\vee}$ for the perturbation problem given by 
\[
\mathcal{F}_{\epsilon}(u)\phi = \mathcal{F}_0(u)\phi + \int_{\Omega} \epsilon_1 \Delta u \Delta \phi 
+ \int_{\Omega} \epsilon_2 N(u)\phi=0.
\]
 
Notice that regular enough solutions of the weak equation above correspond to 
strong solutions (with $\beta=1$) of the problem 
\[
 F_{\epsilon}(u)= \epsilon_1 \Delta^2  u + \beta \Delta u + \alpha u(1-u)  + \epsilon_2 N(u)=0,
\]
with \textit{natural boundary conditions}.

Indeed, starting with
\[
\int_{\Omega} \epsilon_1 \left(\Delta^2 u \right) \phi + \int_{\Omega} \beta (\Delta u)\phi + \int_{\Omega} (\alpha u(1-u) + \epsilon_2 N(u))\phi=0,
\]
and applying Green's formula we have:
\[
\int_{\Omega} -\epsilon_1 \nabla(\Delta u) \cdot \nabla \phi + \oint_{\partial \Omega} \epsilon_1 (\nabla(\Delta u) \cdot n)\phi 
\]
\[
 - \int_{\Omega} \beta \nabla u \cdot \nabla \phi + \oint_{\partial \Omega} \beta (\nabla u \cdot n) \phi
+ \int_{\Omega}  ( \alpha u(1-u) + \epsilon_2 N(u) )\phi=0.
\]
Assuming that the boundary integrals vanish, we apply Green's formula once more and now have:
\[
\int_{\Omega} \epsilon_1 \Delta u \Delta \phi - \oint_{\partial \Omega} \epsilon_1 (\nabla \phi \cdot n) \Delta u 
- \int_{\Omega} \beta \nabla u \cdot \nabla \phi
+ \int_{\Omega} (\alpha u (1-u) + \epsilon_2 N(u))\phi=0.
\]
Noting that the boundary integral vanish, we obtain the weak equations 
\[
\int_{\Omega} \epsilon_1 \Delta u \Delta \phi 
- \int_{\Omega} \beta \nabla u \cdot \nabla \phi
+ \int_{\Omega} (\alpha u(1-u) + \epsilon_2 N(u))\phi=0,
\]
i.e
\[
\int_{\Omega} \epsilon_1 \Delta u \Delta \phi 
- \int_{\Omega} \beta \nabla u \cdot \nabla \phi
+ \int_{\Omega} \alpha u \phi= \int_{\Omega} (\alpha u^2 - \epsilon_2 N(u))\phi.
\]

The main purpose of presenting the simple derivation above is to explicitly state the meaning of \textit{natural boundary conditions} for the problem in consideration. 

In this last form, one easily identify the perturbation problem  
from Fisher's equation, $u_t = F_{\epsilon}(u)$, where

\[
F_{\epsilon}(u)=\epsilon_1 \Delta^2 u + \beta \Delta u + \alpha u(1-u) + \epsilon_2 N(u) = F_0(u) + \epsilon_1 \Delta^2 u + \epsilon_2 N(u).  
\]
We will choose $N(u)=-0.5|\nabla u|^2$ for our computations (and $\beta=1$), and $\epsilon_1$ will be a small negative parameter. In this way, for $\beta=0$ and  $\alpha=0$ (and with Dirichlet boundary conditions instead) we recover the Kuramoto-Sivashinsky model. On the other hand, with $\epsilon_1=0$ and $\epsilon_2=0$ we obtain again Fisher's equation. 

\begin{remark}
The computations in the Matlab scripts are formulated as
$u_t= -\alpha \Delta^2 u -\beta \Delta u + \mu u(1-u) -\delta 0.5|\nabla u|^2,$
with $\alpha$ small and positive, $\beta$ negative of absolute value close to 1, $\mu$ close to the parameters used for Fisher's equation, and $\delta$ small and positive. 
\end{remark}

Because the weak form of the equation contains the Laplacian (instead of just the gradient)
we  use $C^1$ Argyris elements which offer high convergence rate. 
We refer to \cite{MR1930132} for the mathematical theory of the Argyris elements
and to  \cite{MR2528520} for a useful discussion of numerical the implementation.  

We we recall that
the Argyris elements are fifth order polynomials in two space variable constructed as follows.
Define the operators
$L_1 = {id}$, $L_2 = \partial_{10}$,  $L_3 = \partial_{01}$, 
$L_4 = \partial_{20}$, $L_5 = \partial_{11}$ and $L_6 = \partial_{02}$.
For an element $[n_1, n_2, n_3, m_1, m_2, m_3]$ with nodes 
$n_1$, $n_2$ and $n_3$ and midpoints $m_1$, $m_2$ and $m_3$,
the nodal basis $\phi_k^{n_i}$ are defined by 
\[
L_{\ell} (\phi_k^{n_i}(n_j)) = \delta_{ij} \delta_{\ell k},
\]
\[
\frac{\partial}{\partial n} \phi_k^{n_i}(m_j) = 0,
\]
and the basis associated to the midpoints by 
\[
\frac{\partial}{\partial n} \phi^{m_i}(m_j) = \delta_{ij},
\]
\[
L_\ell \phi^{m_i}(n_j) = 0,
\]
where $1 \leq i, j \leq 3$ and $1 \leq k, \ell \leq 6$. 

These are $21$ constraints for each $\phi$ which uniquely defines a fifth 
order polynomials of the form
\[
\phi(x,y) = \sum_{0 \leq i + j \leq 5} c_{ij} x^i y^j.
\]
We solve a $21 \times 21$ linear system $Ac=b$ for the coefficients $c_{ij}$ for each of the 21 basis associated with an element.
In practice, we only do this for a reference triangle and transfer these basis to an arbitrary element 
using the method of Dominguez and Sayas \cite{MR2528520}.

In the notation presented earlier, $S_i=\{L_k: 1\leq k \leq 6\}$ for $i=1,2,3$ and  
$S_i=\{\frac{\partial}{\partial n}\}$ for $i=4,5,6$.  After some indexing and renaming 
we let $S= \bigcup\limits_i S_i = \{L_k\}$ and
\[
\phi_i=\frac{det(A_i)}{det(A)}
\]
for $1\leq i\leq 21$.
The global representation of $u$ becomes:
\[
u = \sum_{k = 1}^6 \sum_{\mbox{\tiny all } n_i} c_k^{n_i} \phi_k^{n_i}
+ \sum_{\mbox{\tiny all } m_i} c^{m_i} \phi^{m_i}.
\]


This interpolation is indexed in some convenient way: $u=\sum\limits_{j=1}^{{nb}}c_j\phi_j$ with ${nb}=6{nv}+{ned}$ where ${nv}$ is the number of vertices and ${ned}$ 
is the number of edges in the triangulation.

After computing an equilibrium solution and eigendata $\lambda$ and $\xi$ as in the previous examples, 
we proceed to solve Equation \eqref{eq:invEq} in the case of Morse index one, interpreted in 
${H^2}^{\vee}$ as
\[
F(P(\theta))=\lambda \theta \frac{\partial P}{\partial \theta}.
\]

\begin{figure}[!t]
\begin{center}
    \begin{minipage}[r]{0.45\textwidth}
	\includegraphics[width=\textwidth]{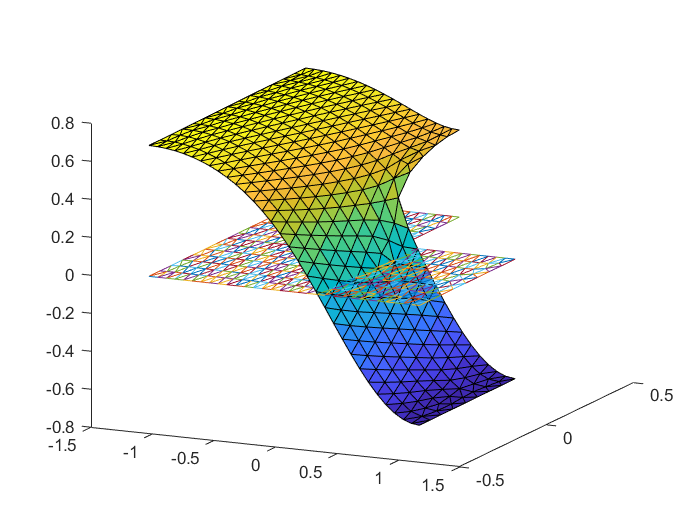}
	\end{minipage}
	\begin{minipage}[l]{0.45\textwidth}
		\includegraphics[width=\textwidth]{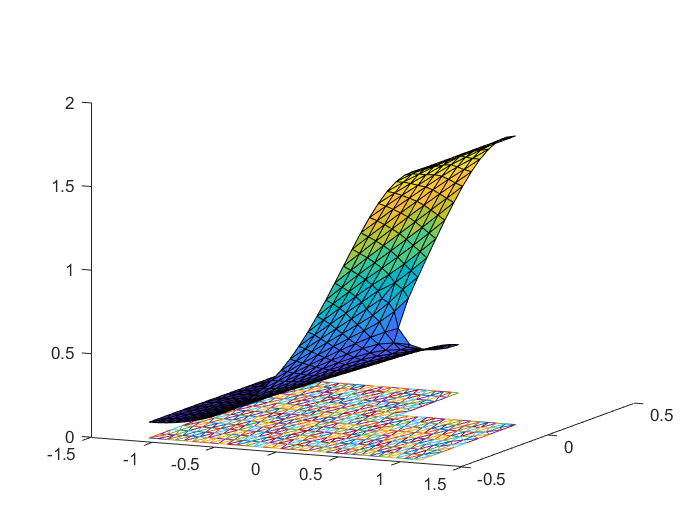}
	\end{minipage}
    \end{center}
    \caption{FKS equation with
    $\epsilon_1=-10^{-2}$, $\beta=1$, $\alpha=2.61$, and $\epsilon_2=10^{-3}$
    Left:  Equilibrium solution, ne=705.
    Right: Eigenfunction $\xi$ with $\lambda=3.48$. }
\end{figure}

First comparing powers and then solving for 
\[
p_n=\sum\limits_{j=1}^{{nb}}c_j^{(n)}\phi_j,
\]
leads to 
\[
\langle \epsilon_1 \Delta^2 p_n +\beta \Delta p_n 
+\alpha(1-2p_0)p_n 
- \epsilon_2 \left(\frac{\partial p_0}{\partial x}\frac{\partial p_n}{\partial x} 
+ \frac{\partial p_0}{\partial y}\frac{\partial p_n}{\partial y} \right)
-\lambda p_n, \phi \rangle = 
\langle s_n, \phi \rangle
\]
where 
\[
s_n=\frac{\epsilon_2}{2} \left(\sum_{k=1}^{n-1}\frac{\partial p_k}{\partial x}\frac{\partial p_{n-k}}{\partial x} 
+ \frac{\partial p_k}{\partial y}\frac{\partial p_{n-k}}{\partial y} \right)
-\alpha \sum_{k=1}^{n-1} p_k p_{n-k},
\]
and so the projected weak formulation of the homological equation is of the form
\[
\left(D \mathcal{F}^{h}(c^{(0)})-\lambda n \int_{\Omega} \phi_j \phi_i \right) c^{(n)} =
\left(\int_{\Omega} s_n \phi_i \right).
\]

\begin{figure}[!t]	
	\begin{center}
		\begin{minipage}[l]{0.45\textwidth}
			\includegraphics[width=\textwidth]{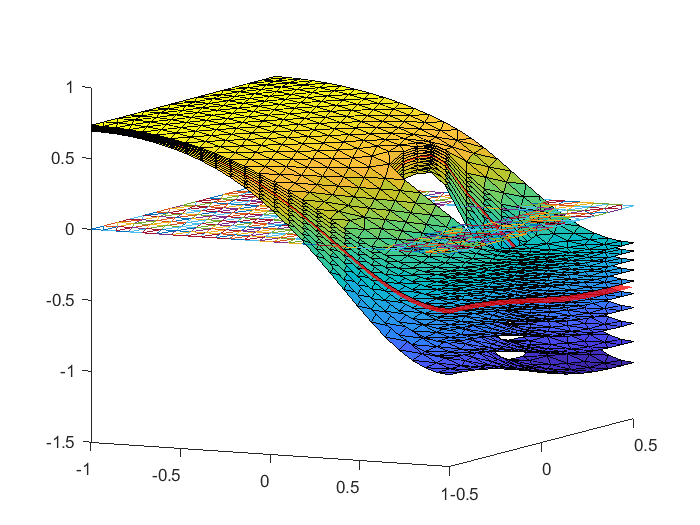}
		\end{minipage}
		\begin{minipage}[r]{0.45\textwidth}
			\includegraphics[width=\textwidth]{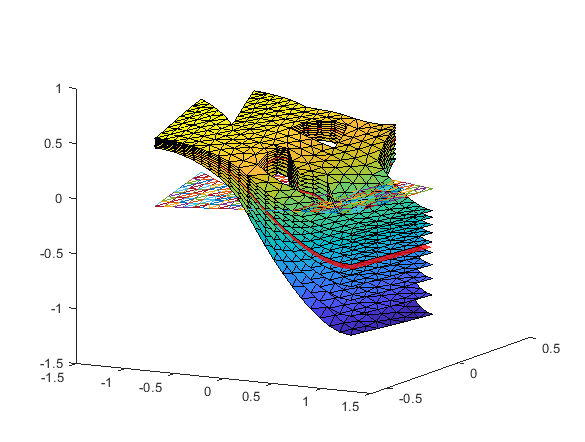}
		\end{minipage}		
	\end{center}
	\caption{Unstable manifolds for the FKS equation posed on 
		non-convex domains with holes. 
		Left: $\epsilon_1=-10^{-3}$, $\beta=1$, $\alpha=3$, and $\epsilon_2=10^{-4}$, 10 points on the 1d manifold, $N=30$. 
		$L^2$ error on the invariance equation 2.09e-07.
		Right:
		$\epsilon_1=-10^{-2}$, $\beta=1$, $\alpha=3$, and $\epsilon_2=10^{-3}$, 10 points on the 1d manifold, $N=30$. 
		$L^2$ error on the invariance equation 1.45e-06.
		} \label{fig:KS_manifolds_differentDomains}
\end{figure}

\begin{figure}[!t]
\begin{center}
    \begin{minipage}[r]{0.45\textwidth}
	\includegraphics[width=\textwidth]{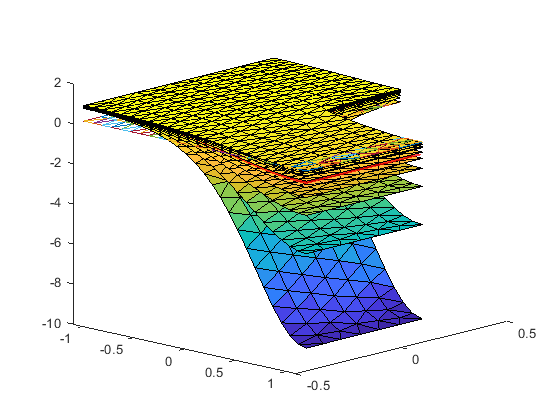}
	\end{minipage}
	\begin{minipage}[l]{0.45\textwidth}
		\includegraphics[width=\textwidth]{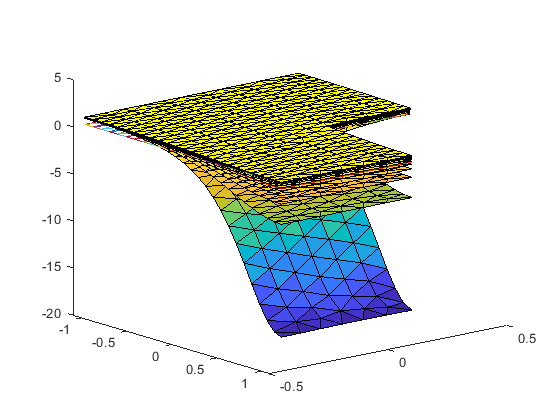}
	\end{minipage}
    \end{center}
    \caption{Unstable manifold for the FKS equation on the   $\mathbb{L}$  domain with eigenvector scaled by 0.8 and parameters
    	$\epsilon_1=-10^{-2}$, $\beta=1$, $\alpha=2.61$, and $\epsilon_2=10^{-3}$
    Left:  10 points on the 1d manifold, $N=10$,
Invariance equation error $L^2$ norm 0.012.
    Right:  10 points on the 1d manifold, $N=120$,
Invariance equation error $L^2$ norm 1.55e-05}
\label{fig:KS_domain_comparison}
\end{figure}

In the Figures \ref{fig:KS_manifolds_differentDomains}, we show the manifolds computed over two additional irregular domains. 
In Figure \ref{fig:KS_domain_comparison} the manifold is approximated to order 10 and 120, 
using the same scaling of the eigenvector. The error improves significantly by increasing the order of the approximation. Equivalently, if we set a tolerance level for the error in our computations, the local manifold obtained for order 10 is significantly smaller. 

%
%
\subsection{A-posteriori error estimation} \label{sec:implementation}

In this section we define a-posteriori error indicators for the parameterization 
method and illustrate their use in the examples from above.  
For the first indicator, consider the $L^2$ norm of the defect associated 
with the invariance equation.  That is, for the $N$-th order 
parameterization 
\[
P^N(\theta, x,y) = \sum_{n = 0}^N p_n(x,y) \theta^n,
\]
of a 1D unstable manifold, define the defect function
\[
E_{1, N}(\theta, x, y):= F(P_N(\theta,x,y))- \lambda \theta \frac{\partial }{\partial \theta} P_N(\theta,x,y),
\] 
for $\theta \in (-1,1)$ and $(x,y) \in \Omega$, 
and the $L^2$ indicator 
\[
\epsilon_{N,1} = ave_{|\theta| \leq 1}\left\| E_{1, N}(\theta) \right\|_{L^2(\Omega)} 
\]

Note that $\epsilon_{N,1} = 0$ for an exact solution.

Similarly we define, for the parameterization 
\[
P^N(\theta_1, \theta_2, x,y) = \sum_{m + n = 0}^N p_{mn}(x,y) \theta_1^m \theta_2^n, 
\] 
of a two dimensional unstable manifold, the defect  function 
\[
E_{N,2}(\theta_1, \theta_2, x, y) = 
F(P_N(\theta_1, \theta_2,x,y)) - 
\lambda_1 \theta_1 \frac{\partial }{\partial \theta_1} P_N(\theta_1, \theta_2,x,y) - 
\lambda_2 \theta_2 \frac{\partial }{\partial \theta_2} P_N(\theta_1, \theta_2,x,y),
\]
and the indicator 
\[
\epsilon_{N,2} = 
ave_{|\theta_1|, |\theta_2| \leq 1}\left\| E_{2, N}(\theta_1, \theta_2) \right\|_{L^2(\Omega)}.
\]

In practice these indicators are approximates by computing the $L^2(\Omega)$ norms and average for
a finite number of parameter points.

Another class of indicators is obtained by considering the dynamical conjugacy error
discussed in Equation \eqref{eq:conj}.  That is, with fixed $T > 0$ define 
the dynamical defect 
\[
\mbox{conjError}(T)_{N,1}(\theta, x,y) = 
P^N(e^{\lambda T} \theta, x,y) - \Phi(P^N(\theta, x,y), T),
\] 
$\theta \in (-1,1)$ and  $(x,y) \in \Omega$, for the 1D manifold and 
\[
\mbox{conjError}(T)_{N,2}(\theta_1, \theta_2, x,y) = 
P^N(e^{\lambda_1 T} \theta_1, e^{\lambda_2 T} \theta_2, x,y) 
- \Phi(P^N(\theta_1,\theta_2, x,y), T), 
\]
$\theta_1, \theta_2 \in [-1,1]\times [-1,1]$ and  $(x,y) \in \Omega$
for the 2D manifold.
Then we have the indicators 
\[
\varepsilon_{N,1} = \sup_{|\theta| \in [-1,1]} \left\|
\mbox{conjError}(T)_{N,1}(\theta)
\right\|_{L^2(\Omega)},
\]
and 
\[
\varepsilon_{N,2} = \sup_{|\theta_1|, |\theta_2| \leq 1} \left\|
\mbox{conjError}(T)_{N,2}(\theta_1, \theta_2)
\right\|_{L^2(\Omega)}
\]
Note that the calculation of these indicators depends on the (fairly arbitrary)
choice of $T$, and more over requires 
numerical approximation of the flow map $\Phi(P(\theta,x,y),t)$, which in 
turn requires implementation of a numerical integration scheme for the 
parabolic PDE.
Then the computation of the $\epsilon$-indicators is in general much simpler 
than the $\varepsilon$-indicators.  For this reason, we much prefer the former 
in the present work.  Nevertheless, the latter can be very valuable for 
debugging purposes, and we always check the conjugacy errors before 
claiming with confidence that we have working codes.

Tables \ref{table:1} - \ref{table:6} below
report the results of a number of defect 
calculations for the manifold computations of the previous section.
We observe that in general the defect decreases as the number of elements 
increases (and hence the mesh size decreases) and tends to improve as the order 
$N$ of the approximation increases. 
It should also be stressed that using finite elements of higher order in a given problem
seems to have a dramatic effect on the error.  This is illustrated
in Table \ref{table:6}, which compares the defect for the 1D manifold in 
the Fisher equation using piecewise linear versus Argyris elements. 
While the piecewise linear elements proved approximately 6 figures of 
accuracy on the L-shaped domain, 
using the higher order elements we obtain defects on the order of 
machine precision.  
The later are considerably more difficult to implement,
but offers significant advantages, and are especially 
encouraging for potential future 
applications in computer assisted proofs.

\begin{table}
\begin{center}
	\begin{tabular}{|c| c| c| c|} 
		\hline
		${ne}$ & Fisher 1d manifold $\mathbb{L}$ & Fisher 2d manifold $\mathbb{L}$\\ 		
		\hline\hline
		515 & 4.922499e-07 & 5.960955e-07
		  \\ 
		\hline
		984 & 1.448931e-07 & 4.655299e-08 \\
		\hline
		1963 & 3.294838e-08 & 4.020046e-08\\
		\hline		
	\end{tabular}
	\end{center}
\caption{
Table: $L^2$ norms of the error in the Invariance Equation for 1d and 2d unstable 
manifolds in the Fisher model over the $\mathbb{L}$ domain: 
$\alpha =2.7$, $\alpha=9$ respectively. }
\label{table:1}
\end{table}

\begin{table}
\begin{center}
	\begin{tabular}{|c| c| c| c|} 
		\hline
		${ne}$ & FR 1d manifold $\mathbb{L}$ & FR 
		2d manifold $\mathbb{L}$ \\ 
		\hline\hline
		515 & 1.804745e-07 & 1.994842e-05  \\ 
		\hline
		984 & 4.655299e-08 & 5.777222e-06 \\
		\hline
		1963 & 1.189424e-08 & 2.417198e-06
		 \\
		\hline		
	\end{tabular}
\end{center}
\caption{
Table: $L^2$ norms of the error in the Invariance Equation for 1d and 2d manifolds in the Fisher model with exponential nonlinearity over the $\mathbb{L}$ domain: $\alpha=-4.7$, $\alpha=-4.41$ respectively. 
} \label{table:2}
\end{table}

\begin{table}
\begin{center}
	\begin{tabular}{|c| c| c|} 
		\hline
		${ne}$ & FKS 1d manifold $\mathbb{L}$ \\
		\hline\hline
		100 & 8.395986e-06
		  \\ 
		\hline
		200 & 4.306036e-06
		\\
		\hline
		423 & 2.208767e-06 \\
		\hline		
	\end{tabular}
\end{center}
\caption{
Table: $L^2$ norms of the error in the Invariance Equation for the 1d unstable manifold over the $\mathbb{L}$ domain: $\epsilon_1=-10^{-2}$, $\beta=1$, $\alpha=2.61$, and $\epsilon_2=10^{-3}$. 
} \label{table:3}
\end{table}

\bigskip

\begin{table}
\begin{center}
	\begin{tabular}{|c| c| c|} 
		\hline
		${ne}$ & FKS 1d manifold Door \\
		\hline\hline
		123 & 7.355615e-07
		 \\ 
		\hline
		260 & 3.895756e-07\\
		\hline
		522 & 2.086379e-07  \\
		\hline
	\end{tabular}
\end{center}
\caption{
Table: $L^2$ norms of the error in the Invariance Equation for the 1d unstable manifold over he door domain: $\epsilon_1=-10^{-3}$, $\beta=1$, $\alpha=3$, and $\epsilon_2=10^{-4}$. }
\label{table:4}
\end{table}

\bigskip

\begin{table}
\begin{center}
	\begin{tabular}{|c| c| c|} 
		\hline
		${ne}$ & FKS 1d manifold Polygon \\
		\hline\hline
		97 & 4.993837e-06  \\ 
		\hline
		193 & 2.897142e-06\\
		\hline
		412 & 1.447564e-06  \\
		\hline
	\end{tabular}
\end{center}
\caption{
Table: $L^2$ norms of the error in the Invariance Equation for the 1d unstable manifold over the polygon with holes: $\epsilon_1=-10^{-2}$, $\beta=1$, $\alpha=3$, and $\epsilon_2=10^{-3}$. 
}\label{table:5}
\end{table}

\bigskip

\begin{table}
\begin{center}
	\begin{tabular}{|c| c| c| c|} 
		\hline
		${ne}$ & Fisher 1d manifold piecewise linear & Fisher 1d manifold Argyris \\ 
		\hline\hline
		423 & 
		6.208993e-07 & 5.777960e-16  \\ 
		\hline
	\end{tabular}
\end{center}
\caption{
Table: $L^2$ norms of the error in the Invariance Equation for 1 dimensional manifolds in the Fisher model over the $\mathbb{L}$ domain with piecewise linear and Argyris basis: $\alpha=2.7$. }
\label{table:6}
\end{table}

\section{Conclusions} \label{sec:conclusions}
We have combined the parameterization method with finite element 
analysis to obtain a new approximation method for 
unstable manifolds of equilibrium solutions for parabolic PDEs.
The method is applied to several PDEs defined  on planar polygonal domains
and is implemented for number of example problems
with both polynomial and non-polynomial nonlinearities, for unstable manifolds of dimension 
one and two, for a number of non-convex and non-simply connected domains, and for problems 
involving both Laplacian and bi-harmonic Laplacian diffusion operators.  
The method is easy to implement for computing the approximation to arbitrary order:
the same code that computes the second order approximation will compute the 
approximation to order $50$ -- this is just a matter of changing a loop variable.
The method is amenable to a-posteriori analysis of errors and we employ these 
indicators to show that our calculations are accurate far from the equilibrium solution.  

Interesting future projects would be to apply the method to problems with other boundary 
conditions such as Dirichlet or Robin, to apply it to problems formulated 
on spatial domains of dimension 3 or more, to extend the method for the computation of 
unstable manifolds attached to periodic solutions of parabolic PDEs, or 
to extend the method to study invariant manifolds attached to equilibrium 
or periodic solutions of systems of parabolic PDEs.  

Finally we mention that there is a thriving literature on mathematically rigorous 
computer assisted proof for elliptic PDEs based on finite element analysis.  
See for example the works of 
\cite{MR944817,MR1284278,MR1203360,MR1131108,MR1151049,MR2378291,MR1354660,MR2431603,MR2109916,MR2019251,MR2406180,MR1354655,MR3203775,MR2492179,MR3196951,MR3611667}
for validated numerical methods for solving nonlinear elliptic PDE
(equilibrium solutions of parabolic PDEs) and their associated 
eigenvalue/eigenfunction problems.   We refer also the references just cited for more
complete review of this literature.  
From the point of view of the present discussion the important point is this:  that 
the present work reduces the problem of computing jets of unstable manifolds
to the problem of solving elliptic boundary value problems -- and moreover that 
a number of authors have developed powerful methods of computer assisted proof
for solving such problems.
A very interesting line of future research would be to combine the results of the 
present work validated numerical methods for elliptic BVPs.

\section{Acknowledgements}
The authors would like to thank Rafael de la Llave, Michael Plum,
and Allan Hungria for helpful 
discussions as this work evolved.  
J.G. and J.D.M.J. were partially supported by 
the Sloan Foundation Grant FIDDS-17.  J.D.M.J. was partially 
supported by the National Science Foundation grant  DMS - 1813501.

\section{Data Availability Statement}
The data that support the findings of this study are available on request from the corresponding author J.D.M.J.

\bibliographystyle{siamplain}      
\bibliography{papers,references}   

\end{document}